\documentclass[a4paper,12pt]{amsart}

\usepackage{geometry}

\usepackage{amsmath}
\usepackage{amsthm}
\usepackage{amsfonts}
\usepackage{amssymb}
\usepackage{mathtools}
\let\amslrcorner\lrcorner
\usepackage{mathtools}
\usepackage{slashed}
\usepackage[dvipsnames]{xcolor}
\usepackage{tikz-cd}
\usepackage{graphicx}
\usepackage{wrapfig}
\usepackage{accents}
\usepackage{hyperref}
\usepackage{enumitem}
\usepackage{bookmark}
\usepackage{mathrsfs}
\usepackage[nameinlink,noabbrev]{cleveref}
\usepackage{dutchcal}
\usepackage[super]{nth}
\usepackage{stmaryrd}
\usepackage{witharrows}
\usepackage{verbatim}
\usepackage{multirow}
\usepackage{pifont}
\usepackage[dvipsnames]{xcolor}

\hypersetup{colorlinks, linkcolor=blue, citecolor=blue}

\usepackage{mathabx}
\let\lrcorner\amslrcorner

\newcommand\restr[2]{{
  \left.\kern-\nulldelimiterspace 
  #1 
  \vphantom{\big|} 
  \right|_{#2} 
  }}

  \newcommand{\RN}[1]{%
	  \textup{\uppercase\expandafter{\romannumeral#1}}%
  }

\makeatletter
\renewcommand*\env@matrix[1][\arraystretch]{%
  \edef\arraystretch{#1}%
  \hskip -\arraycolsep
  \let\@ifnextchar\new@ifnextchar
  \array{*\c@MaxMatrixCols c}}
\makeatother

\newcommand{\V}[1]{\textup{#1}}

\newcommand{\dir}{\slashed{\partial}}

\newcommand{\fall}{\hspace{.15cm} \forall \hspace{.08cm}}

\newcommand{\sth}{\hspace{.15cm} | \hspace{.15cm}}

\newcommand{\Sreg}{\mathcal{R}_{\V{reg}}}

\DeclareMathOperator{\im}{Im}

\DeclareRobustCommand{\frcshape}{\fontfamily{frc}\selectfont}
\DeclareTextFontCommand{\textfrc}{\frcshape}

\DeclarePairedDelimiter\abs{\lvert}{\rvert}%
\DeclarePairedDelimiter\norm{\lVert}{\rVert}%

\theoremstyle{plain}
\newtheorem{theorem}{Theorem}[section]
\newtheorem*{theorem*}{Theorem}
\newtheorem*{maintheorem*}{Main Theorem}

\newtheorem{proposition}[theorem]{Proposition}
\newtheorem*{proposition*}{Proposition}
\newtheorem{lemma}[theorem]{Lemma}

\theoremstyle{definition}

\newtheorem{definition}[theorem]{Definition}
\newtheorem*{definition*}{Definition}

\newtheorem{observation}[theorem]{Observation}
\newtheorem*{observation*}{Observation}

\theoremstyle{remark}
\newtheorem{remark}{Remark}[section]

\newtheorem*{remark*}{Remark}

\newtheoremstyle{dotlessP}{}{}{\color{violet}}{}{\color{violet}\bfseries}{}{ }{}
\theoremstyle{dotlessP}

\def\XXint#1#2#3{{\setbox0=\hbox{$#1{#2#3}{\int}$ }
\vcenter{\hbox{$#2#3$ }}\kern-.6\wd0}}

\numberwithin{equation}{section}

\title{Transversality for Perturbed Special Lagrangian Submanifolds}
\author{Emily Autumn Windes}

\begin{document}
\maketitle

\begin{abstract}
	We prove a transversality theorem for the moduli space of perturbed special Lagrangian submanifolds in a 6-manifold equipped with a generalization of a Calabi-Yau structure. These perturbed special Lagrangian submanifolds arise as solutions to an infinite-dimensional Lagrange multipliers problem which is part of a proposal for counting special Lagrangians outlined by Donaldson and Segal in \cite{dose09}. More specifically, we prove that this moduli space is generically a set of isolated points. 
\end{abstract}

\section{Introduction}

The prospect of extending gauge-theoretic ideas, originally developed for dimensions 2, 3 and 4, to manifolds with special holonomy in dimensions 6, 7 and 8 was explored by Donaldson and Thomas in \cite{doth96} and by Donaldson and Segal in \cite{dose09}. In particular, manifolds with special holonomy come equipped with calibrations, whose corresponding calibrated submanifolds are of great interest to both physicists and mathematicians. In certain respects, these calibrated submanifolds are analogous to \( J\)-holomorphic curves and exhibit many connections to other gauge-theoretic objects such as instantons and monopoles. Donaldson, Thomas, and Segal put forth proposals for how one might count various calibrated submanifolds in the hope of developing new invariants for manifolds with special holonomy. 

There are several technical reasons such a program is not straightforward. This direction was considered by Joyce \cite{joyc16}, Doan and Walpuski \cite{dowa17}, and others. This paper is concerned specifically with the situation in dimension 6. 

It was discovered by Hitchin \cite{hitc00, hitc01} that metrics with special holonomy in dimensions 6 and 7 are deeply connected to the existence and properties of differential forms whose pointwise \( \V{GL}(n,\mathbb{R})\)-orbit is open. Such forms are called \textit{stable}. One very interesting feature of dimension 6 is the fact that an \( \V{SU}(3)\)-structure on a 6-manifold is equivalent to a choice of a stable 3-form \( \rho\) and a stable 4-form \( \tau\) satisfying certain algebraic conditions (see \Cref{def:SU(3)-structure}). Furthermore, a Calabi-Yau structure on a 6-manifold is an \( \V{SU}(3)\)-structure \( (\rho,\tau)\) where the 3-form \( \rho\) and the 4-form \( \tau\) along with their Hitchin duals (see \cref{stable-forms-dim-6}) are \textit{closed}. In the Calabi-Yau setting, the 3-form \( \rho\) is the real part of the holomorphic volume form and its Hitchin dual \( \hat{\rho}\) is the imaginary part. 

In a 6-dimensional Calabi-Yau manifold, both \( \rho\) and \( \hat{\rho}\) are calibrations, whose calibrated submanifolds are called \textit{special Lagrangian}. In \cite{dose09}, Donaldson and Segal pointed out that special Lagrangian submanifolds can be characterized as solutions to a certain Lagrange multipliers problem defined purely in terms of these stable forms. We briefly summarize this Lagrange multipliers problem here. See \cref{lmp} for a more rigorous exposition. 

Suppose that \( (M,\rho,\tau)\) is a manifold equipped with a pair of closed, stable forms \( (\rho,\tau) \in \Omega^3(M) \times \Omega^4(M)\). Any Calabi-Yau manifold satisfies this condition, but as we will see, we will also consider more general settings. Next, fix a 3-submanifold \( L_{0} \subset M\). Let \( L\) be any nearby 3-submanifold representing the same homology class as \( L_0\) and let \( Q\) be a cobordism connecting \( L_0\) and \( L\). Then, roughly speaking, one can define a functional \( f_{\tau}\) on the space of submanifolds by integrating the 4-form \( \tau\) over \( Q\)
\begin{equation*}
	f_{\tau}(L) = \int_{Q} \tau.
\end{equation*}
Of course, this functional is not well-defined since it depends on the choice of cobordism \( Q\). However, it is well-defined on the covering space of the space of embeddings of a 3-manifold into \( M\) as explained in \cref{lmp}.

Next, let \( C_{\rho}\) be the set of 3-submanifolds in \( M\), all diffeomorphic to \( L_0\), representing the same homology class and satisfying \( \restr{\rho}{L} = 0\). Note that in the Calabi-Yau case, a special Lagrangian submanifold which is calibrated by the 3-form \( \hat{\rho}\) always satisfies this condition. The Donaldson--Segal Lagrange-multipliers problem is to find the critical points of \( f_{\tau}\) restricted to the set \( C_{\rho}\). 

In a Calabi-Yau manifold, special Lagrangian submanifolds are critical points of the Lagrange functional arising from this Lagrange multipliers problem. This fact raises the question of whether it is possible to construct a Floer theory for special Lagrangian submanifolds of a 6-dimensional Calabi-Yau manifold. 

However, a result of McLean \cite{mcle98} says that for any special Lagrangian submanifold \( L\), the moduli space of nearby special Lagrangians is a smooth manifold with dimension equal to the first Betti number of \( L\). Thus, with the goal of defining a Floer theory in mind, we want to perturb the \( \V{SU}(3)\)-structure underlying the Calabi-Yau structure on our 6-manifold so that the critical points of the relevant functional are isolated. 

The fact that 7 is 6+1 gives us a natural space of perturbations. This is because any 6-manifold \( M\) with an \( \V{SU}(3)\)-structure can be embedded as a hypersurface in a 7-dimensional cylinder \( \mathbb{R} \times M\) equipped with a \( G_2\)-structure arising in a canonical way from the \( \V{SU}(3)\)-structure. Specifically, if \( (\rho,\tau) \in \Omega^3(M) \times \Omega^4(M)\) defines an \( \V{SU}(3)\)-structure on \( M\), then 
\begin{equation*}
	\psi = dt \wedge \rho + \tau
\end{equation*}
is a stable 4-form on \( \mathbb{R} \times M\) corresponding to a \( G_2\)-structure on \( \mathbb{R} \times M\). If the \( \V{SU}(3)\)-structure on \( M\) is a Calabi-Yau structure, the metric associated to the 4-form \( \psi\) on \( \mathbb{R} \times M\) will have holonomy \( \V{SU}(3)\) contained in \( G_2\). On the other hand, if \( (M,\rho,\tau)\) is a 6-manifold equipped with a pair of stable forms \( (\rho,\tau) \in \Omega^3(M) \times \Omega^4(M)\), the condition that this pair gives rise to a \( G_2\)-structure on \( \mathbb{R} \times M\) is a \textit{weaker} condition than the requirement that it defines an \( \V{SU}(3)\)-structure on \( M\). Let 
\begin{equation*}
	\mathcal{R}_{G_2} = \left\{ (\rho,\tau) \in \Omega^3(M) \times \Omega^4(M) \sth d\rho = 0, d\tau = 0, \text{and }\psi = dt \wedge \rho + \tau \text{ is stable on } \mathbb{R} \times M \right\}
\end{equation*}
be the set of \textit{closed} \( G_2\)-\textit{pairs} on a 6-manifold \( M\).

If we start with a Calabi-Yau 6-manifold \( M\), we can perturb the underlying \( \V{SU}(3)\)-structure \( (\rho,\tau)\) in such a way that it is no longer an \( \V{SU}(3)\)-structure but nonetheless gives rise to a \( G_2\)-structure on \( \mathbb{R} \times M\). The reason we require \( \rho\) and \( \tau\) to be closed is to ensure that the deformation theory of the \( \psi\)-associative submanifolds in \( \mathbb{R} \times M\) is governed by a \textit{self-adjoint} operator which is not necessarily the case unless \( \psi\) is closed. Since the condition that a form be stable is an \textit{open condition}, any small perturbation of the of the pair \( (\rho,\tau)\) by closed forms remains in the space of closed \( G_2\)-pairs since such a perturbation corresponds to a small perturbation of the 4-form \( \psi\). As long as \( \rho\) and \( \tau\), and hence the 4-form \( \psi\), are closed, the Lagrange-multipliers problem described above can still be defined. The solutions to this more general Lagrange-multipliers problem generalize the notion of a special Lagrangian. They are no longer calibrated submanifolds since we no longer require that the Hitchin duals of \( \rho\) and \( \tau\) are closed but nonetheless retain several useful features typically enjoyed by calibrated submanifolds.

The Euler-Lagrange equations corresponding to the Donaldson--Segal Lagrange-multipliers problem take the following form.

\begin{proposition}
	Let \( (M,\rho,\tau)\) be a 6-manifold equipped with a \( G_2\)-pair \( (\rho,\tau)\) where both \( \rho\) and \( \tau\) are closed. A 3-submanifold \( L \subset M\) is a solution to the Lagrange multipliers problem if and only if there exists \( \lambda \in C^{\infty}(L)\) such that 
	\begin{eqnarray*}
		\tau_N + d\lambda \wedge \rho_N &=& 0 \\ 
		\restr{\rho}{L} &=& 0.
	\end{eqnarray*}
\end{proposition}

These equations first appeared in \cite{dose09} without much explanation. The details of their derivation can be found in \cref{perturbed-SL-equations}. See \Cref{def:aN} for an explanation of the notation. Here, they will be referred to as the \textit{perturbed special Lagrangian equations} (perturbed SL equations). The function \( \lambda\) is the \textit{Lagrange-multiplier} and plays a similar role to the Lagrange-multipliers found in calculus. When \( M\) is a Calabi-Yau manifold with holomorphic volume form \( \Omega\) and symplectic structure \( \omega\), the 3-form \( \rho = \V{Re}(\Omega)\) and the 4-form \( \tau = \frac{1}{2}\omega^2\). In this case, the solutions to the perturbed SL equations are precisely the special Lagrangian submanifolds together with constant functions \( \lambda\). When \( M\) is not necessarily Calabi-Yau, a submanifold \( L\) solving the perturbed SL equations for some function \( \lambda\) is called a \textit{perturbed special Lagrangian submanifold}.

The perturbed SL equations can also be considered separately from the Lagrange-multipliers set-up. In this case, it is not necessary to require that \( \rho\) and \( \tau\) be closed. One could still hope to construct a numerical invariant for 6-manifolds equipped with a \( G_2\)-pair in this way, although this is not the focus of this paper. 

The key connection between solutions to the perturbed SL equations and the 7-dimensional setting is as follows. Suppose that \( (M,\rho,\tau)\) is a 6-manifold with a \( G_2\)-pair and that \( (\lambda,L)\) is a solution to the perturbed SL equations. Then the \textit{graph} of the function \( \lambda\) over \( L \subset M\) is an \textit{associative submanifold} (see \Cref{def:associative}) of \( \mathbb{R} \times M\). For this reason, we call the pair \( (\lambda,L)\) a \textit{graphical associative}. See \Cref{lem:associative-graph} for proof of this fact. As a consequence, the perturbed special Lagrangian equations are \textit{elliptic} since the deformation theory of associative submanifolds in a \( G_2\)-manifold is governed by an elliptic operator \cite{mcle98}. 

The \textit{moduli space} \( \mathcal{M}(A,P;(\rho,\tau))\) is defined to be the set of all submanifolds \( L\) in \( M\), diffeomorphic to a particular 3-manifold \( P\) and representing a fixed homology class \( A \in H_3(M;\mathbb{Z})\), which satisfy the perturbed SL equations for some \( \lambda\).

The definition of a \( G_2\)-pair is flexible enough to prove transversality for \( \mathcal{M}(A,P;(
\rho,\tau))\). 
More precisely, let \( P\) be a closed, oriented 3-manifold and \( M\) a closed, oriented 6-manifold equipped with a closed \( G_2\)-pair \( (\rho,\tau) \in \Omega^3(M) \times \Omega^4(M)\). We prove

\begin{theorem}[transversality]
	Fix a homology class \( A \in H_3(M;\mathbb{Z})\) and a closed 3-manifold \( P\). There is a residual subset \( \mathcal{R}_{\V{reg}}\) of \( \mathcal{R}_{G_2}\) such that the moduli space \( \mathcal{M}(A,P;(\rho,\tau))\) is a collection of isolated points whenever \( (\rho,\tau) \in \mathcal{R}_{\V{reg}}\).
	\label{thm:main-theorem}
\end{theorem}

When the \( G_2\)-pair \( (\rho,\tau)\) is of class \( C^{\ell}\), the proof of this theorem is a consequence of the implicit function theorem for Banach spaces. In order to extend the result to smooth \( G_2\)-pairs, we must apply the so-called Taubes trick which is a standard method from symplectic geometry (see for example chapter 3 of \cite{mcsa04}). In order to utilize the Taubes trick, we must prove an elliptic regularity theorem (\Cref{thm:regularity}) and a compactness theorem (\Cref{thm:compactness}) for associative submanifolds with bounded second fundamental form and bounded volume. This is non-trivial because the solutions to the perturbed SL equations are not necessarily minimal submanifolds. 

If the elements of the moduli space \( \mathcal{M}(A,P;(\rho,\tau))\) are to be counted, then we must determine if or when the moduli space is compact. Towards this end, we we adapt the definition of a \textit{tamed} \( G_2\)-structure (first introduced in \cite{dose09}) to define a \textit{tamed} \( G_2\)-pair on \( M\). This is analogous to the notion of a tamed almost complex structure first introduced by Gromov in order to control the energy of \( J\)-holomorphic curves. In our case, a \( G_2\)-pair \( (\rho,\tau)\) is \textit{tamed} by a second pair \( (\rho',\omega') \in \Omega^3(M) \times \Omega^2(M)\) where both \( \rho'\) and \( \omega'\) are closed and stable. See \cref{tamed-structures} for details. We have 

\begin{proposition}[volume bounds]
	Let \( (M,\rho',\omega')\) be a closed, 6-manifold equipped with taming forms \( (\rho',\omega') \in \Omega^3(M) \times \Omega^2(M)\). Suppose that \( (\rho,\tau)\) is a \( (\rho',\omega')\)-tame \( G_2\)-pair. Then every element in \( \mathcal{M}(A,P;(\rho,\tau))\) has a topological volume bound depending only on the homology class of \( \rho'\).
\end{proposition}

More details about tamed \( G_2\)-structures can be found in \cite{joyc16}. Although in general we cannot expect that the second fundamental form of a perturbed SL submanifold to be bounded, we nonetheless expect the volume bound in the previous proposition to give us a compactification of the moduli space using rectifiable currents. However, such a compactification must necessarily contain singular objects which are not entirely understood. Such singular objects were studied by Joyce in \cite{joyc03}. We hope to explore these topics along with other topics related to developing a Floer theory for special Lagrangian submanifolds in future papers.

\subsection{Organization}
The paper is organized as follows. \Cref{stable-forms} includes a review of results about stable forms in dimensions 6 and 7 as well as a brief discussion about the different contexts in which this work might be applied. \Cref{ta-g2} contains information about \( G_2\)-structures on cylinders of the form \( \mathbb{R} \times M\) where \( M\) is a 6-manifold equipped with a pair of stable forms. The notion of a \( G_2\)-pair is defined and its relationship to \( \V{SU}(3)\)-structures is discussed. In \cref{tamed-structures}, we adapt the definition of a tamed \( G_2\)-structure to define a tamed \( G_2\)-pair. \Cref{lmp} contains a detailed description of the Lagrange multipliers problem first appearing in \cite{dose09}. The perturbed SL equations for this Lagrange multipliers problem are defined. In \cref{ellipticity}, we show that solutions to the perturbed SL equations correspond to associative submanifolds in a 7-dimensional cylinder. We then use this fact to prove that the perturbed SL equations are elliptic. In \cref{volume-bounds} we derive volume bounds for perturbed SL submanifolds in a 6-manifold with a tamed \( G_2\)-pair. \Cref{elliptic-regularity} contains regularity and compactness results for associative submanifolds. Finally, in \cref{transversality} we apply all of the above to prove a transversality result for the moduli space of perturbed special Lagrangian submanifolds. 

\subsection{Acknowledgments} This paper would not have been possible without the continuous support, guidance, and expertise of my supervisors Aleksander Doan and Boris Botvinnik. I also thank Thomas Walpuski, Jason Lotay, Robert Bryant, Lorenzo Foscolo, and Costante Bellettini for several enlightening discussions on topics related to this paper. Lastly, I thank Jesse Madnick and Oliver Edtmair for helpful suggestions and feedback on rough drafts of this paper. 

\section{Stable forms} \label{stable-forms}

\begin{definition}
	A form \( \phi \in \Lambda^{p}(\mathbb{R}^n)^{*}\) is called \textit{stable} if its \(\V{GL}(n,\mathbb{R})\)-orbit in \(\Lambda^{p}(\mathbb{R}^n)^*\) is open. 
	\label{def:stable_forms}
\end{definition}

Stable forms and their connections to metrics with special Holonomy were first studied in \cite{hitc01}. Precisely which manifolds admit stable forms and in which degrees was worked out in \cite{hpva07}. There it was proved that stable 3-forms only exist in dimensions 6, 7, and 8. In the present paper, we will be mainly interested in stable forms in dimension 6, but will rely heavily on the connection to special geometry in 7 dimensions. It should be noted however that although there exist stable 3-forms in dimension 8, metrics on 8-manifolds with \( \V{Spin}(7)\) holonomy do not arise from them in the same way that metrics with special holonomy arise from stable 3-forms on 7-manifolds.

If a manifold \(M\) admits a stable \(p\)-form \(\alpha\), then it has a \(G\)-structure where \(G\) is the pointwise stabilizer of \( \alpha\). Let \( \alpha\) represent a stable 2-form in any dimension, or a stable 3-form in dimension 6, 7, or 8. In both these cases, the stabilizer group preserves a volume form denoted by \( \V{vol}(\alpha)\). It was discovered by Hitchin that metrics with special holonomy in both 6- and 7-dimensions have variational characterisations in terms of the functional which inputs a stable form and outputs the total volume of the manifold with respect to a choice of an invariant volume form. Next, we recall features of stable forms in 6 and 7 dimensions in more detail. 

\subsection{Dimension 6}\label{stable-forms-dim-6}

In dimension 6, stable 2-forms, 3-forms, and 4-forms are possible. Their stabilizers and conventional volume forms are as follows:

\begin{enumerate}
	\item A stable 2-form \(\omega\) is a non-degenerate 2-form. The group \( \V{Sp}(6,\mathbb{R})\) can be defined to be the set of linear transformations of \( \mathbb{R}^6\) that preserve a non-degenerate 2-form, and is a real, non-compact, connected, simple Lie group. The associated volume form is 
		\[ \V{vol}(\omega) = \frac{1}{6} \omega^3\] 
		which is also known as the \textit{Liouville volume form}.
	\item A stable 4-form \(\tau\) also has stabilizer \(\V{Sp}(6,\mathbb{R}) \). To define \( \V{vol}(\tau)\), use the isomorphism \( I:\Lambda^{4}(\mathbb{R}^6)^* \cong \Lambda^{2} \mathbb{R}^6 \otimes \Lambda^{6}(\mathbb{R}^{6})^{*}\). Then \( I(\tau)^{6} \in \Lambda^{6}\mathbb{R}^{6} \otimes \left( \Lambda^{6} (\mathbb{R}^6)^{*} \right)^{3} \cong \left( \Lambda^{6}(\mathbb{R}^{6})^{*} \right)^{2}\). Define 
		\[ \V{vol}(\tau) = \left( I(\tau)^{6} \right)^{1/2}.\]
	\item A stable 3-form \(\rho\) has stabilizer \(\V{SL}(3,\mathbb{C}) \) (in which case it's called \textit{positive}) or stabilizer \( \V{SL}(3,\mathbb{R}) \times \V{SL}(3,\mathbb{R})\) (in which case it's called \textit{negative}). Fix a 3-form \( \rho \in \Lambda^{3}(\mathbb{R}^6)^*\). Define a map \( K_{\rho}(v) = v \lrcorner \rho \wedge \rho \in \Lambda^{5}(\mathbb{R}^6)^* \cong \mathbb{R}^{6} \otimes \Lambda^6(\mathbb{R}^6)^*\). Then the positive forms are those forms with \( \V{tr}(K)^2 < 0.\) In this case, we define

		\[ \V{vol}(\rho) = \abs*{\sqrt{-\V{tr}(K)^2}} \in \Lambda^6 (\mathbb{R}^6)^* .\]
\end{enumerate}

In this paper, we will only be concerned with stable 3-forms whose orbit is \( \V{SL}(3,\mathbb{C})\). From now on, a whenever a 3-form is referred to as stable, we actually mean positive. For more information about how the volume forms are defined, see the appendix of \cite{hitc01}. In the same paper, Hitchin uses the homogeneous behavior of the map from stable forms to volume forms to define a \textit{Hitchin dual} which for a stable form \( \alpha\) will be denoted by \( \hat{\alpha}\). 

Another special feature of dimension 6 is that a positive 3-form determines an almost complex structure, a fact that was explained in \cite{brya06}. The Hitchin duals for stable 2- and 4-forms, and for positive 3-forms on \( \mathbb{R}^6\) are described below explicitly. 

\begin{enumerate}
	\item For a stable 2-form \(\omega\), \(\hat{\omega} = \frac{1}{2} \omega^2\). 
	\item For a positive 3-form \(\rho\), \(\hat{\rho}\) is the unique 3-form \(\hat{\rho}\) such that \(\rho + i \hat{\rho}\) is a nowhere vanishing complex volume form with respect to the complex structure determined by \( \rho\). 
	\item For a stable 4-form \(\tau\), \(\hat{\tau}\) is the unique nondegenerate 2-form satisfying \(\tau = \frac{1}{2} \hat{\tau}^2\).
\end{enumerate}

In fact, we can use the concept of stable forms to \textit{define} \(\V{SU(3)}\) structures on 6-manifolds in much the same way that \( G_2\) structures on \( 7\)-manifolds are often said to be defined by a stable 3-form.

\begin{definition}
	An \(\V{SU(3)}\)-structure on a 6-manifold \(M\) is a pair of differential forms \( (\omega, \rho)\) such that 
	\begin{enumerate}
		\item \(\omega\) is a stable 2-form 
		\item \(\rho \) is a positive 3-form 
		\item the following algebraic conditions are satisfied:
			\begin{equation*}
				\omega \wedge \rho = 0, \qquad \qquad \frac{1}{6} \omega^3 = \frac{1}{4} \rho \wedge \hat{\rho}.
			\end{equation*}
	\end{enumerate}

	In this case, the stabilizer of the pair \((\omega, \rho)\) will be exactly \(\V{SU(3)}\). An \( \V{SU}(3)\)-structure can be equivalently defined in terms of a stable 4-form and a stable (positive) 3-form, where the 2-form appearing in the last condition is the Hitchin dual of the 4-form. Throughout the paper, \( \omega\) and \( \omega'\) will always denote a 2-form, \( \tau\) will always denote a 4-form, and \( \rho\) and \(\rho'\) will always denote a 3-form. 
	\label{def:SU(3)-structure}
\end{definition}

Various conditions may be placed on the stable forms, most of which have been studied in some detail. Suppose that \( (M,\rho,\tau)\) is a manifold with an \( \V{SU}(3)\)-structure. 

\begin{enumerate}
	\item If both \( \rho\) and \( \hat{\rho}\) are closed then \( M\) is a \textit{complex threefold with a trivial canonical bundle}. Manifolds such as these are not necessarily Calabi-Yau manifolds because they are not necessarily K\"ahler. These manifolds were studied in \cite{tosa14} where they are referred to as non-K\"ahler Calabi-Yau manifolds. There are many simple examples such as \( S^1 \times S^3\). 
	\item If both \( \rho\) and \( \tau\) are closed, then \( M\) is said to have a \textit{half-flat} \( \V{SU}(3)\)-structure. Half-flat structures are important to the study of hypersurfaces in \( \mathbb{R}^{7}\) and were first studied by Calabi. They have since showed up again in physics and geometry. A slightly more restricted sub-class of six manifolds with half-flat \( \V{SU}(3)\)-structures are the \textit{nearly K\"ahler} manifolds whose Riemannian cones have holonomy equal to \( G_2\). See \cite{fosc17} for more details about nearly K\"ahler manifolds and \cite{masa13} for an explanation of half-flat \( \V{SU}(3)\)-structures. 
	\item If \( \hat{\tau}\) is closed, then \( M\) is a \textit{symplectic manifold with a compatible almost-complex structure}. Note that \( \hat{\tau}\) being closed also implies that \( \tau\) is closed since \( \tau = \frac{1}{2} \hat{\tau}^2\). 
	\item If \( \tau, \hat{\tau}, \rho,\) and \( \hat{\rho}\) are all closed, then \( M\) is a \textit{Calabi-Yau} manifold. 
\end{enumerate}

Note that we could also drop the requirement that \( (\rho,\tau)\) forms an \( \V{SU}(3)\) structure and simply study manifolds that carry a pair of stable forms. In this paper, we will usually consider pairs of stable forms that determine a co-closed \( G_2\) structure in \( 6 +1\) dimensions. More details can be found in the following sections. 

\begin{remark} 
	Recall that the existence of an almost complex structure implies the existence of a non-degenerate (i.e. stable) 2-form, since one can be constructed out of an almost complex structure and any Riemannian metric. This is true in any dimension. Therefore in dimension 6, the existence of a stable 3-form implies the existence of a stable 2-form since a stable 3-form defines an almost complex structure.
\end{remark}

Next we review the 7-dimensional situation. 

\subsection{Dimension 7} Fix an identification \( \mathbb{R}^{7} \cong \V{Im}\mathbb{O}\) of \( \mathbb{R}^7\) with the imaginary octonions. Then the multiplication on \( \mathbb{O}\) endows \( \mathbb{R}^{7}\) with

\begin{enumerate}
	\item an inner product \( g\) 
		\begin{equation*}
			g(u,v) := -\V{Re}(uv)
		\end{equation*}
	\item a cross-product

		\begin{equation*}
			u \times v := \V{Im}(uv)
		\end{equation*}
	\item a 3-form
		\begin{equation*}
			\varphi(u,v,w):= g(u \times v, w)
		\end{equation*}
	\item an \textit{associator} \( [\cdot,\cdot,\cdot]: \Lambda^{3}\mathbb{R}^{7} \rightarrow \mathbb{R}^{7}\)
		\begin{equation}
			[u,v,w] := (u \times v) \times w + \langle v,w\rangle u - \langle u,w\rangle v
			\label{eq:associator}
		\end{equation}
	\item and a 4-form 
		\begin{equation*}
			\psi(u,v,w,z) := g\left( [u,v,w],z \right).
		\end{equation*}
\end{enumerate}

One can always choose coordinates \( x_1,\dots,x_7\) on \( \mathbb{R}^{7}\) so that the 3-form \( \varphi\) above can be written as
\begin{equation*}
	\varphi_0 = dx_{123} + dx_{145} + dx_{167} + dx_{246} - dx_{257} - dx_{347} - dx_{356}.
	\label{eq:standard-G2-structure}
\end{equation*}
The notation \( dx_{i_1,\dots,i_k}\) denotes the wedge product \( dx_{i_1} \wedge \dots \wedge dx_{i_k}\). The following relations always hold
\begin{eqnarray*}
	\left( u \lrcorner \varphi \right) \wedge \left( v \lrcorner \varphi \right) \wedge \varphi &=& 6g(u,v)\V{vol}_{g} \\
	*_{g} \varphi &=& \psi.
\end{eqnarray*}

\begin{definition}
	The subgroup of \(\V{GL}(7,\mathbb{R})\) which fixes \(\varphi_0\) is the compact, connected, simple Lie group \(G_2\) 
	\begin{equation*}
		G_2 = \left\{ A \in \V{GL}(7,\mathbb{R}) \sth A^* (\varphi_0) = \varphi_0 \right\}.
		\label{eq:G2-group}
	\end{equation*}
	\label{def:G2-group}
\end{definition}

The group \(G_2\) also preserves the standard metric and orientation when acting on \(\mathbb{R}^7\). In particular, \(G_2\) also fixes the 4-form \(\psi_0 = *\varphi_0\), where the Hodge star is taken with respect to the metric determined by \( \varphi\). These forms are stable with respect to \Cref{def:stable_forms} since the Lie group \(G_2\) is 14-dimensional and \(\text{dim}\V{GL}(7,\mathbb{R}) - \text{dim}\Lambda^3(\mathbb{R}^7) = 49 - 35 = 14\), implying that their \( \V{GL}(7,\mathbb{R})\)-orbit is open. From this perspective, the Hitchin dual of \( \varphi\) is \( \psi\). There also exist stable 3-forms on \( \mathbb{R}^{7}\) with stabilizer split \( G_2\). These are analogous to the negative stable 3-forms in dimension 6, and will not be considered in this paper. Throughout, a stable 3-form on a 7-dimensional manifold will always be one whose pointwise stabilizer is the compact real Lie group \( G_2\).

It should be noted that the stabilizer of the 4-form \( \psi = *\varphi\) is actually \( \pm G_2 = G_2 \cup -\V{Id}G_2\). Therefore a stable 4-form on a 7-dimensional vector space determines an inner-product but \textit{not} an orientation. 

\begin{definition}
	A 7-dimensional manifold \(X\) equipped with a global, stable 3-form \(\varphi\) has a \(G_2\)-structure and is called a \textit{\(G_2\)-manifold}. If \(\varphi\) is closed, then \(X\) is called a manifold with a \textit{closed} \(G_2\)-structure. If \(*\varphi = \psi\) is closed, then \(X\) is called a manifold with a \textit{co-closed} \(G_2\)-structure. If \( \varphi\) is both closed and co-closed, then \( X\) is called a manifold with a \textit{torsion-free} \( G_2\)-structure.
\end{definition}

Manifolds with a \( G_2\) structure \( \varphi\) are automatically equipped with a metric \( g_{\varphi}\), a cross-product \( \times_{\varphi}\), an associator \( \left[ \cdot,\cdot,\cdot \right]_{\varphi}\), and a stable 4-form \( \psi = *_{\varphi} \varphi\). Equivalently, we may refer to a 4-form \( \psi\) on an oriented 7-manifold as a \( G_2\) structure. The following result is due to Fernandez and Gray. 

\begin{lemma}[See Lemma 11.5 of \cite{sala89}]
	Let \( X\) be a 7-dimensional manifold equipped with a \( G_2\)-structure \( \varphi\). Let \( g\) be the associated metric. Then the following are equivalent.
	\begin{enumerate}
		\item \( \V{Hol}(g) \subseteq G_2\) 
		\item \( d\varphi = d^*\varphi= 0\)
	\end{enumerate}
\end{lemma}

If \( g\) is a metric with holonomy contained in \( G_2\), \( \varphi\) is the corresponding stable 3-form, and \( \nabla\) is the Levi-Civita connection of \( g\) then \( \nabla \varphi = 0\). The form \( \nabla \varphi\) is called the torsion of \( \varphi\), so this lemma justifies the above terminology. For more details about \(G_2\) structures and \(G_2\) manifolds, see for example \cite{brya05}. \\

\section{Translation-invariant \texorpdfstring{\( G_2\)}{okaaaa}-structures} \label{ta-g2}

In this section, we investigate how \( G_2\)-structures on a cylinder of the form \( \mathbb{R} \times M\) where \( M\) is a 6-manifold relate to pairs of stable forms on \( M\). This discussion will also allow us to set some terminology that will be used in the following sections.

The Lie group \( \V{SU}(3)\) is a subgroup of the Lie group \( G_2\). Indeed, if one chooses a non-zero vector \( v_{0} \in \mathbb{R}^{7}\) then the subgroup of \( G_2\) that preserves \( v_{0}\) is isomorphic to \( \V{SU}(3)\). Therefore there are interesting relationships between manifolds with \( G_2\)-structures and manifolds with \( \V{SU}(3)\)-structures. Throughout this section, let \( M\) be a 6-manifold, \( X = \mathbb{R} \times M\), and \( (\rho, \omega) \in \Omega^{3}(M) \times \Omega^{2}(M)\). 

\begin{definition} 

	A pair of stable forms on a 6-manifold \( M\) that define a \( G_2\)-structure on the cylinder \( \mathbb{R} \times M\) will be called a \( G_2\)-\textit{pair}. These pairs can consist of a stable 3-form and a stable 2-form \textit{or} a stable 3-form and a stable 4-form. In the sections that follow, the Greek letter \( \omega\) will \textit{always} denote a 2-form, the Greek letter \( \rho\) will \textit{always} denote a 3-form, and the Greek letter \( \tau\) will \textit{always} denote a 4-form. Set 
	\begin{equation*}
		\mathcal{R}_{G_2} = \left\{ (\rho,\omega) \in \Omega^3(M) \times \Omega^2(M) \sth \varphi = \rho + dt \wedge \omega \in \Omega^{3}_{\V{stable}}(X) \right\}
	\end{equation*}
	\textit{or}
	\begin{equation*}
		\mathcal{R}_{G_2} = \left\{ (\rho,\tau) \in \Omega^3(M) \times \Omega^4(M) \sth \psi = dt \wedge \rho + \tau \in \Omega^{4}_{\V{stable}}(X) \right\}.
	\end{equation*}
\label{def:G2-pair}
\end{definition}

\begin{remark}
	A pair of stable forms \( (\rho,\omega) \in \Omega^3(M) \times \Omega^2(M)\) is a \( G_2\) pair if and only if the pair \( (\hat{\rho},\tau)\) of their Hitchin duals is a \( G_2\) pair. One can see this by showing that if \( \psi = dt \wedge \hat{\rho} + \tau\) is stable then the metric associated to \( \varphi = dt \wedge \omega + \rho\) must be non-degenerate and vice-versa.
\end{remark}

\begin{definition} Similarly, we will let \( \mathcal{R}_{\V{SU}(3)}\) denote the set of \( \V{SU}(3)\)-structures on \( M\). That is 
	\begin{equation*}
		\mathcal{R}_{\V{SU}(3)} = \left\{ (\rho,\omega) \in \Omega^3(M) \times \Omega^2(M) \sth (\rho,\omega)\ \text{is an SU}(3)\text{-structure}  \right\}
	\end{equation*}
	\textit{or}
	\begin{equation*}
		\mathcal{R}_{\V{SU}(3)} = \left\{ (\rho,\tau) \in \Omega^3(M) \times \Omega^4(M) \sth (\rho,\tau)\ \text{is an SU}(3)\text{-structure}  \right\}.
	\end{equation*}
	\label{def:SU3-pair}
\end{definition}

\begin{remark}
	In the case of \( \V{SU}(3)\)-structures, a pair \( (\rho,\tau)\) is exactly equivalent to a pair \( (\rho,\omega)\) where \( \tau\) is the Hitchin dual of \( \omega\) since in this case, the metric defined on the cylinder is the product metric and therefore the same whether or not we choose the 4-form or the 2-form.
\end{remark}

	The set \( \mathcal{R}_{\V{SU}(3)}\) is a subset of \( \mathcal{R}_{G_2}\) and they are not the same as the following standard lemmas illustrate. 

\begin{lemma}
	The pair \( (\rho,\omega) \in \Omega^{3}(M) \times \Omega^{2}(M)\) is an \( \V{SU}(3)\)-structure on \( M\) with metric \( g_M\) if and only if \( \varphi = \rho + dt \wedge \omega\) is a \( G_2\)-structure that induces the product metric \( g_{\varphi} = dt^2 + g_{M}\) on \( X = \mathbb{R} \times M\). 
	\label{lem:tiG2-1}
\end{lemma}

\begin{lemma}
	The pair \( (\rho,\omega)\) defines a Calabi-Yau structure on \( M\) if and only if \( \varphi = \rho + dt \wedge \omega\) is a torsion-free \( G_2\)-structure on \( X = \mathbb{R} \times M\). In this case, \( \V{Hol}(g_{\varphi}) \subseteq \V{SU}(3)\).
	\label{lem:tiG2-2}
\end{lemma}

The proofs of Lemmas \ref{lem:tiG2-1} and \ref{lem:tiG2-2} follow from the proof of Proposition 11.1.1 in \cite{joyc00}.

\begin{remark}
	There are cases where both \( \rho \in \Omega^3(M)\) and \( \omega \in \Omega^2(M)\) are stable, but \( \varphi = \rho + dt \wedge \omega\) is not. For example, let \( M = \mathbb{R}^{6}\) with coordinates \( x_1, \dots , x_6\). Set 
\begin{eqnarray*}
	\rho &=& dx_{135}+dx_{632} + dx_{254} + dx_{416} \\
	\omega &=& dx_{63} + dx_{25} + dx_{41}.
\end{eqnarray*}
One can check that \( \rho\) and \( \omega\) are stable. However, 
\begin{equation*}
	\omega \wedge \rho = dx_{63254} + dx_{25416} + dx_{41632} \neq 0
\end{equation*}
so \( (\rho,\omega)\) is \textit{not} an \( \V{SU}(3)\)-structure. Next, let \( x_0\) denote an extra coordinate spanning \( \mathbb{R}\). Set 
\begin{eqnarray*}
	\varphi &=& \rho + dx_0 \wedge \omega \\
	&=& dx_{135}+dx_{632} + dx_{254} + dx_{416} + dx_{063} + dx_{025} + dx_{041}. 
\end{eqnarray*}
Then one can compute that 
\begin{equation*}
	(\partial_{x_1} \lrcorner \varphi) \wedge \left( \partial_{x_1} \lrcorner \varphi \right) \wedge \varphi = 0
\end{equation*}
and thus \( g_{\varphi}\) is not a definite form, so \( \varphi\) is not stable.
\end{remark}

\begin{remark}
	There are cases where \( (\rho,\omega) \in \Omega^3(M) \times \Omega^2(M)\) is \textit{not} an \( \V{SU}(3)\)-structure but \( \varphi = \rho + dt \wedge \omega\) \textit{does} define a \( G_2\)-structure. In this case, the metric \( g_{\varphi}\) does \textit{not} define a product metric by \Cref{lem:tiG2-1}. To see this, let \( \rho = dx_{135} + dx_{632} + dx_{254} + dx_{416}\) and \( \omega = dx_{12} + dx_{34} + dx_{56}\). Then 
\begin{equation*}
	\varphi= dx_{135}+dx_{632} + dx_{254} + dx_{416} + dx_{012} + dx_{034} + dx_{056}. 
\end{equation*}
Since being stable is an open condition, small modifications of \( \varphi_0\) are still stable. For example, 
\begin{equation*}
	\varphi' = dx_{135}+dx_{632} + dx_{254} + dx_{416} + dx_{012} + dx_{034} + dx_{056} + Kdx_{123}
\end{equation*}
where \( K\) is a constant is still stable for small enough \( K\). Note that this is equivalent to changing \( \rho\) out for \( \rho + Kdx^{123}\). Note also that 
\begin{equation*}
	\left( \rho + Kdx_{123} \right) \wedge \omega = Kdx_{12356} \neq 0
\end{equation*}
so after this modification, \( (\rho,\omega)\) do \textit{not} form an \( \V{SU}(3)\)-structure. Furthermore, a computation shows that 
\begin{equation*}
	(\partial_{x_3} \lrcorner \varphi')\wedge \left( \partial_{x_0} \lrcorner \varphi' \right) \wedge \varphi' = 6g_{\varphi'}(\partial_{x_3},\partial_{x_0})\V{vol}_{g_{\varphi'}} \neq 0
\end{equation*}
but, for any product metric \( g = dx_{0}^{2} + g_{M}\), we have
\begin{equation*}
	g(\partial_{x_3},\partial_{x_0}) = dx_{0}^{2}(0,\partial_{x_0}) + g_{M}(\partial_{x_3},0) = 0.
\end{equation*}
Therefore, \( g_{\varphi'}\) does not define a product metric on \( \mathbb{R} \times M\) for any metric \( g_{M}\) on \( M\). 
\label{ob:tiG2-1}
\end{remark}

\begin{remark}
	Suppose that \( (\rho,\omega)\in \Omega^3(M) \times \Omega^2(M)\) is a \( G_2\) pair, but not necessarily an \( SU(3)\) structure. Then there exists a stable 2-form \( \omega'\) such that \( (\rho,\omega') \) is an \( \V{SU}(3)\)-structure on \( M\). This follows from the fact that \( M = \{0\} \times M \) is a hypersurface in \( \mathbb{R} \times M\). Let \( \times_{\varphi}\)  denote the cross-product associated to \( \varphi\). Then, given any unit normal vector field \( n\) along \( M\), we can define an almost complex structure 
	\begin{equation*}
		J_nv = n \times_{\varphi} x \qquad \qquad \text{for all\ } v \in TM,
	\end{equation*}
	and a stable 3-form given by
	\begin{equation*}
		\rho = \iota^*\varphi
	\end{equation*}
	where \( \iota:M \rightarrow \mathbb{R} \times M\) is the inclusion. Let \( \omega'\) be defined by 
	\begin{equation*}
		\omega'(u,v) = g_{\varphi}(u,-J_nv) \qquad \qquad \text{for all\ } u,v \in TM.
	\end{equation*}
	The pair \( (\rho,\omega')\) defines an \( \V{SU}(3)\)-structure. The \( \V{SU}(3)\)-structures of this type were first studied in \cite{cala58} and more recently by in \cite{cabr04}. We emphasize that, in example given after \Cref{ob:tiG2-1}, the vector field \( \partial_{x_0}\) is \textit{not} normal to \( M\) with respect to the metric \( g_{\varphi'}\). 
\label{lem:tiG2-3}
\end{remark}

\begin{lemma}
	The space \( \mathcal{R}_{\V{SU}(3)}\) is a deformation-retract of \( \mathcal{R}_{G_2}\).
\end{lemma}

\begin{proof}
	Clearly \( \mathcal{R}_{\V{SU}(3)} \subseteq \mathcal{R}_{G_2}\). Let \( (\rho,\omega) \in \mathcal{R}_{G_2}\). Let \( t\) denote the \( \mathbb{R}\)-coordinate on \( \mathbb{R} \times M\) and set \( \varphi = \rho + dt \wedge \omega\) as usual. Choose a unit vector field \( n\) along \( \left\{ 0 \right\} \times M\) that is normal to \( \left\{ 0 \right\} \times M\) with respect to the metric \( g_{\varphi}\). Note that even if \( g_{\varphi}\) is not a product metric on \( \mathbb{R} \times M\), the vector field \( \partial_t\) is nonvanishing, and therefore homotopic to \( n\). Let \( \left\{ n_s \right\}_{s=0}^{1} \) be a smooth family of nowhere-vanishing vector fields satisfying \( n_0 = \partial_t\) and \( n_1 = n\). If \( \varphi\) \textit{does} define the product metric, simply set \( n = \partial_t = n_s\) for all \( s \in [0,1]\). 
	Then for each \( s\) define
	\begin{equation*}
		\rho = \iota^*\varphi \qquad J_s(\cdot) = n_s \times_{\varphi} \cdot \qquad \omega_s(\cdot,\cdot) = g_{\varphi}\left( \cdot, -J_{s}(\cdot)\right).
	\end{equation*}
	Since there exists a normal vector field for each pair \( (\rho,\omega) \in \mathcal{R}_{G_2}\), this construction defines a map
	\begin{equation*}
		F:\mathcal{R}_{G_2} \times I \rightarrow \mathcal{R}_{G_2}
	\end{equation*}
	by 
	\begin{equation*}
		F\left( (\rho,\omega),s \right) = (\rho,\omega_s).
	\end{equation*}

	Clearly \( F\left( (\rho,\omega),0 \right) = (\rho,\omega)\). On the other hand, \( F\left( (\rho,\omega),1 \right) \in \mathcal{R}_{\V{SU}(3)}\) (see also proposition 4.1 of \cite{cabr04}). Clearly if \( (\rho,\omega) \in \mathcal{R}_{\V{SU}(3)}\), then \( F\left( (\rho,\tau),1 \right)\) is the identity. Thus, \( F\) is a deformation retract as desired. 
\end{proof} 

The above lemma also applies to \( G_2\)-pairs of the form \( (\rho,\tau)\) where \( \tau\) is a 4-form. \\ 


\section{Tamed structures} \label{tamed-structures}

In this section we develop the notion of a tamed \( G_2\)-pair. The idea of a \textit{tamed structure} was first introduced by Gromov in \cite{grom85} in order to control the energy of \( J\)-holomorphic curves in a symplectic manifold. Let \(M\) be a smooth manifold equipped with a closed, stable 2-form \(\omega\) and an almost complex structure \(J\). The almost complex structure is called \(\omega\)\textit{-tame} if it satisfies \(\omega(v, Jv) > 0\) for every nonzero tangent vector \(v \in TM\). Let \( (\Sigma, j)\) be a Riemann surface. Then a map \( u:\Sigma \rightarrow M\) satisfying \( J \circ du = du \circ j\) is called a \( J\)-holomorphic curve. Such maps are hugely important to symplectic geometry \cite{mcsa04}.

The energy of a \( J\)-holomorphic curve \( u\) is defined to be \(E(u) = \frac{1}{2} \int_{\Sigma} \abs*{du}_{J}^2 \V{vol}_{\Sigma}\), where the norm \( \abs{du}^{2}_{J}\) depends on the choice of \( J\). When \(J\) is \(\omega\)-tame, one can prove the energy identity
\begin{equation*}
	E(u) = \int_{\Sigma} u^* \omega
\end{equation*}
which is a topological invariant depending only on the homology class represented by \(u\). See Lemma 2.2.1 of \cite{mcsa04} for more details about the symplectic case. This result is essential to proving compactness results for moduli spaces of \( J\)-holomorphic curves.  

In \cite{dose09} a similar notion of \textit{taming} and \textit{tamed} forms was introduced for manifolds with special holonomy in 6, 7 and 8 dimensions. This definition is most easily stated for the 7-dimensional case. We will adapt this definition to 6-manifolds equipped with a \( G_2\)-pair. In our situation, the role of \( J\)-holomorphic curves is played by \textit{\( \psi\)-associative} submanifolds of a 7-manifold. 

Let \( X = \mathbb{R}^7\) and fix an identification of \( \mathbb{R}^7\) with the imaginary octonions. Then the associator (from \cref{eq:associator}) vanishes precisely on the 3-dimensional subspaces of \( X\) that are associative, hence the name.  

\begin{definition}
	Suppose that \((X,\varphi,\psi)\) is an oriented 7-manifold with a (not necessarily torsion-free) \( G_2\)-structure. Let \([\cdot, \cdot, \cdot]\) denote the associator corresponding to \(\psi\) and  \(x \in X\). We say that an oriented 3-plane \(V \subset T_x X\) is \textit{\(\psi\)-associative} if \(\restr{[\cdot, \cdot, \cdot]}{V} \equiv 0\) and \(\restr{\varphi}{V} > 0\). Similarly, if \(\iota: P \rightarrow X\) is a 3-submanifold of \(X\), then \(\iota P\) is called \textit{\(\psi\)-associative} if \(\iota^*[\cdot, \cdot, \cdot] \equiv 0\) and \(\iota^* \varphi > 0\). 
	\label{def:associative}
\end{definition}

The following definition allows us to relax the requirement that a \( G_2\)-structure be torsion-free, but retain volume-boundedness of associative submanifolds. For an in-depth explanation of the benefits and drawbacks of this condition as well as several equivalent conditions, see section 2.6 of \cite{joyc16}. 

\begin{definition}[Donaldson--Segal]
	Suppose that \((X,\psi)\) is a 7-manifold equipped with a stable 4-form \( \psi\). Then we say that a closed 3-form \(\varphi'\) \textit{tames} \(\psi\) if for all \(x \in X\), and for all \(\psi\)-associative, oriented 3-planes \(V \subset T_x X\), there exists a positive constant \(K\) such that 
	\begin{equation*}
		\V{vol}_{V} \leq K \restr{\varphi'}{V}
	\end{equation*}
	where \( \V{vol}_{V}\) denotes the volume form on \( V\)  with respect to the metric induced by the one given by the \( G_2\)-structure associated to \( \psi\).
	\label{def:tamed-G2-structure}
\end{definition}

If \( \psi\) is co-closed, then \( \varphi = *\psi\) always tames \( \psi\). Note that the definition of a taming 3-form does not stipulate that it must be stable. However, the taming condition does imply this. 

\begin{lemma}
	Let \( (X,\psi)\) be a 7-manifold equipped with a stable 4-form. Suppose that the 3-form \( \varphi'\) tames \( \psi\). Then \( \varphi'\) is stable. 
	\label{lem:tame-implies-stable}
\end{lemma}

\begin{proof}
	It suffices to prove this for \( X = \mathbb{R}^{7}\). Let \( \psi\) be a stable 4-form on \( \mathbb{R}^7\). Suppose that a 3-form \( \varphi'\) tames \( \psi\), but is \textit{not} a stable form. Note that a stable 3-form \( \varphi\) determines a symmetric, bilinear form on \( \mathbb{R}^7\) via the formula
	\begin{equation*}
		6 g_{\varphi}(u,v)\V{vol}_{{\varphi}} = (u \lrcorner \varphi) \wedge (v \lrcorner \varphi) \wedge \varphi.
	\end{equation*}
	This form is positive-definite precisely when \( \varphi\) is positive, and negative-definite precisely when \( \varphi\) is negative, so if \( \varphi'\) is not stable, then \( g_{\varphi'}\) is neither positive definite nor negative definite. That means there exists a nonzero \( u \in \mathbb{R}^7\) such that 
	\begin{equation*}
		(u \lrcorner \varphi') \wedge (u \lrcorner \varphi') \wedge \varphi' = 0.
	\end{equation*}
	We show that this contradicts the assumption that \( \varphi'\) tames \( \psi\). Note that if \( u \lrcorner \varphi' = 0\) then \( \varphi'\) clearly does not tame \( \psi\) since \( u\) is contained in some \( \psi\)-associative 3-plane. There are two remaining cases.

	\textbf{Case 1.} Suppose that \( (u \lrcorner \varphi') \wedge (u \lrcorner \varphi')= 0\). Let \( \left\{ e_i \right\}_{i = 1}^{7} \) be a \( g_{\varphi'}\)-orthogonal basis for \( \mathbb{R}^7\) with \( e_1 = u\) and let \( \left\{ e^{i} \right\}_{i =1}^{7}\) be the \( g_{\varphi'}\)-dual basis. Then \( (u \lrcorner \varphi') \wedge (u \lrcorner \varphi') = 0 \) implies that \( u \lrcorner \varphi'\) must be indecomposable. That is, \( u \lrcorner \varphi' =A e^{ij}\) for some constant \( A\) with \( i, j \neq 1\). Next, let \( v \in \mathbb{R}^{7}\) be a nonzero vector orthogonal to the 3-plane spanned by \( u, e_i\), and \( e_j\) with respect to the metric \( g_{\varphi'}\). Then \( V = \V{span}\left\{ u,v, u \times_{\psi} v \right\}\), where \( \times_{\psi}\) is the cross-product defined by the 4-form \( \psi\), is \( \psi\)-associative. But \( v \lrcorner u \lrcorner \varphi' = v \lrcorner Ae^{ij} = 0\) since \( v\) is orthogonal to \( e_i\) and \( e_j\). This means that \( \restr{\varphi'}{V} = 0\) so \( \varphi'\) does not tame \( \psi\). 

	\textbf{Case 2.} Suppose that \( (u \lrcorner \varphi') \wedge (u \lrcorner \varphi') \wedge \varphi' = 0\) but \( (u\lrcorner \varphi') \wedge (u \lrcorner \varphi') \neq 0\). Let \( \beta = (u \lrcorner \varphi') \wedge (u \lrcorner \varphi') \). This is impossible because the map sending \( \beta\) to \( \beta \wedge \varphi'\) is an isomorphism (see Proposition 2.2.1 in \cite{kari05}). 

\end{proof}

We can extend \Cref{def:tamed-G2-structure} to 6-manifolds as follows. 

\begin{definition}
	Let \( (\rho,\tau) \in \Omega^{3}(M) \times \Omega^{4}(M)\) be a \( G_2\) pair as in \Cref{def:G2-pair}. We say that the pair \((\rho', \omega') \in \Omega^{3}_{\V{closed}}(M) \times \Omega^{2}_{\V{closed}}(M)\) \textit{tames} \((\rho,\tau)\) if 
	\begin{eqnarray*}
		\psi &=& \tau + dt \wedge \rho \\
		\varphi' &=& \rho' + dt \wedge \omega'
	\end{eqnarray*}
	comprises a \textit{tamed} \( G_2\)-structure as in \Cref{def:tamed-G2-structure}.
	\label{def:tamed-SU(3)-structure-1}
\end{definition}

As in the 7-dimensional case, if \( (\rho',\omega')\) tame a \( G_2\)-pair, then both \( \rho'\) and \( \omega'\) must be stable. Thus they are also a \( G_2\)-pair. The condition that a pair tames another is always an open condition. In fact, the set of \( G_2\) 3-forms which tame a \( G_2\) 4-form is an open, convex cone (see Proposition 2.8 in \cite{joyc16}). For example, in the setting that is of most interest to us, if we start out with a 6-manifold \( (M,\rho,\tau)\) equipped with an \( \V{SU}(3)\)-structure and consider the associated 7-manifold \( X = \mathbb{R} \times M\) with the associated \( G_2\)-structure \( \psi=\tau + dt\wedge \rho\), then not only does \( *\psi\) always tame \( \psi\) but any nearby 3-form  \textit{near} \( *\psi\) also tames \( \psi\). Equivalently, \( G_2\)-pair near \( (\hat{\rho},\omega)\) tames \( (\rho,\tau)\). 

\section{The Lagrange multipliers problem}\label{lmp}

The purpose of this section is to describe the Lagrange multipliers problem that is at the center of this paper. As we will see, the solutions to this Lagrange multipliers problem are special Lagrangian submanifolds in the special case where we start with a Calabi-Yau manifold. These ideas first appeared in \cite{dose09}. The material from the previous section will be used to show that when the relevant structure is tamed, the solutions to the following Lagrange multipliers problem have bounded volume. 

\subsection{Set-up} Much of the notation in what follows is taken from \cite{sawa10} and adapted to this 6-dimensional setting. Suppose that \( M\) is a closed 6-manifold equipped with a \( G_2\)-pair \( (\rho,\tau) \in \Omega^{3}(M) \times \Omega^{4}(M)\). Suppose also that both \( \rho\) and \( \tau\) are closed, but do not require that their Hitchin duals are closed. Fix a closed, oriented, 3-manifold \( P\) and homology class \( A \in H_{3}(M;\mathbb{Z})\). Let 
\begin{equation*}
	\mathcal{F} = \left\{ \iota:P\rightarrow M \sth \iota\ \text{is smooth embedding } , [\iota] = A,\ \iota^{*}\hat{\rho} > 0 \right\}.
\end{equation*}
The tangent space to \( \mathcal{F}\) at \( \iota\) is 
\begin{equation*}
	T_{\iota}\mathcal{F} = \Gamma\left( \iota^*TM \right).
\end{equation*}
Let \( \mathcal{G}\) denote the group of orientation-preserving diffeomorphisms of \( P\) and define 
\begin{equation*}
	\mathcal{S} = \mathcal{F}/\mathcal{G}
\end{equation*}
which can be identified with the space of oriented, 3-dimensional submanifolds of \( M\) diffeomorphic to \( P\) along which \( \hat{\rho}\) is positive. Let \( [\iota]\) denote the equivalence class in \( \mathcal{S}\) of an element \( \iota \in \mathcal{F}\). The tangent space \( T_{[\iota]}\mathcal{S}\) is the quotient 
\begin{equation*}
	T_{[\iota]}\mathcal{S} = \frac{\Gamma(\iota^*TM)}{\left\{ d\iota \circ X \sth X \in \Gamma(TP) \right\}}. 
\end{equation*}

Since \( (\rho,\tau)\) determines a \( G_2\) structure on \( \mathbb{R} \times M\), it also determines a metric on \( \mathbb{R} \times M\) and therefore an induced metric on \( M\). This tangent space therefore may be identified with space of normal vector fields on \( \iota P\) in \( M\). We will let \( N\iota\) denote the normal bundle of \( \iota P\) whenever we want to take this perspective. Next, fix a particular embedding \( \iota_{0}\) of \( P\) into \( M\). Let \( \tilde{\mathcal{F}}\) denote the \textit{universal cover} of \( \mathcal{F}\) based at \( \iota_{0}\). That is, 
\begin{equation*}
	\tilde{\mathcal{F}} = \left\{ \tilde{\iota}: [0,1] \times P \rightarrow M \sth \tilde{\iota}(0,\cdot) = \iota_{0},\ \tilde{\iota}(t,\cdot) = \iota_{t} \in \mathcal{F} \fall t \in [0,1] \right\} / \sim
\end{equation*}
where \( \tilde{\iota} \sim \tilde{\iota}\hspace{.02cm}'\) if \( \tilde{\iota}\) and \( \tilde{\iota}\hspace{.02cm}'\) have the same endpoints and are smoothly homotopic. The group of orientation-preserving diffeomorphisms \( \mathcal{G}\) also has a covering space \( \tilde{\mathcal{G}}\), which is the group of smooth isotopies from \( [0,1]\) to \( \V{Diff}(P)\) starting at the identity. Let 
\begin{equation*}
	\tilde{\mathcal{S}} = \tilde{\mathcal{F}}/\tilde{\mathcal{G}}.
\end{equation*}
Define a functional, \( f_{\tau}: \tilde{\mathcal{F}} \rightarrow \mathbb{R}\) by 
\begin{equation*}
	f_{\tau}(\tilde{\iota}) = \int_{[0,1] \times P} \tilde{\iota}\hspace{.02cm}^*\tau.
\end{equation*}
This functional is well-defined since \( \tau\) is closed. Its derivative \( df_{\tau}\) is a one-form on \( \mathcal{F}\) given by 
\begin{equation*}
	\left( df_{\tau} \right)_{\iota}(n) = \int_{P} \iota^*\left( n\lrcorner \tau \right)
\end{equation*}
for all \( n \in T_{\iota}\mathcal{F}\). Note that both \( f_{\tau}\) and its differential are gauge invariant in the sense \( f_{\tau}\left( \tilde{g}^*\tilde{\iota} \right) = f_{\tau}(\tilde{\iota})\) for any \( \tilde{g} \in \tilde{\mathcal{G}}\). Similarly, \( \left( df_{\tau}\right)_{g^{*}\iota}\left( g^{*}n \right) = \left( df_{\tau} \right)_{\iota}(n)\) for any \( g \in \mathcal{G}\). Also note that if \( n\) and \( n'\) are in the same equivalence class in \( T_{[\iota]}\mathcal{S}\), then \( \left( df_{\tau} \right)_{\iota}(n) = \left( df_{\tau} \right)_{\iota}(n')\) since \( \left( df_{\tau} \right)_{\iota}(v) = 0\) for all \( v \in \Gamma(T\iota P)\). Thus \( f_{\tau}\) descends to a functional on \( \tilde{\mathcal{S}}\) and \( df_{\tau} \) descends to a one-form on \( \mathcal{S}\).

Next, we define the constraint. Let \( c: \mathcal{F} \rightarrow \Omega^{3}(P)\) denote the function given by 
\begin{equation*}
	c(\iota) = \iota^*{\rho}.
\end{equation*}
Similarly, let \( \tilde{c}:\tilde{\mathcal{F}} \rightarrow \Omega^{3}(P)\) denote the function given by 
\begin{equation*}
	\tilde{c}(\tilde{\iota}) = \iota_{1}^{*}(\rho).
\end{equation*}
Then let 
\begin{equation*}
	C_{\rho} = c^{-1}(0) \quad \text{and} \quad \tilde{C}_{\rho} = \tilde{c}^{-1}(0).
\end{equation*}
Note that \( C_{\rho}\) is also gauge invariant in the sense that if \( \iota^*\rho = 0\), then all elements in the equivalence class \( [\iota]\) also satisfy this condition. Let 
\begin{equation*}
	\mathcal{C}_{\rho} = C_{\rho} / \mathcal{G} \quad \text{and} \quad \tilde{\mathcal{C}}_{\rho} = \tilde{C}_{\rho} / \tilde{\mathcal{G}}.
\end{equation*}
Next, we want to show that \( \mathcal{C}_{\rho}\) is a submanifold of \( \mathcal{S}\). First, a lemma. 

\begin{lemma}
	Suppose that \( [\iota] \in \mathcal{C}_{\rho}\) so that \( \iota^{*}\rho = 0\) and \( \iota^*\hat{\rho} > 0\). Then for every \( p \in P\), the map \( N_{\iota(p)}\iota \rightarrow \Lambda^2\left( T_{p}^{*} P \right)\) given by 
	\begin{equation*}
		n \mapsto \iota^*(n \lrcorner \rho)
	\end{equation*}
	is an isomorphism.
	\label{lem:isomorphism}
\end{lemma}

\begin{proof}
	Let \( (\rho,\omega')\) be any \( \V{SU}(3)\)-structure on \( M\) and let \( g\) be the corresponding metric with respect to this \( \V{SU}(3)\)-structure. There is an orthogonal splitting of the k-forms on \( M\) with respect to \( g\). We fix the following notation for this splitting throughout the proof. Let \( N\) denote the normal bundle of \( \iota P\) and \( T\) denote the tangent bundle. Let \( p \in P\) and \( x = \iota(p)\). 
	\begin{eqnarray*}
		\Lambda^2\left( T_x^{*}M \right) &=& \Lambda^2(N^*_x) \oplus \left( N_x^{*} \otimes T_x^{*} \right) \oplus \Lambda^2\left( T_x^{*} \right) \\
		&=& \Lambda^{2,0} \oplus \Lambda^{1,1} \oplus \Lambda^{0,2} \\
		\Lambda^3\left( T_x^* M \right) &=& \Lambda^3(N_x^*) \oplus \left( \Lambda^2\left( N^*_x \right) \otimes T_x^* \right) \oplus \left( N_x^* \otimes \Lambda^2\left( T_x^* \right) \right) \oplus \Lambda^3\left( T_x^* \right) \\
		&=& \Lambda^{3,0} \oplus \Lambda^{2,1} \oplus \Lambda^{1,2} \oplus \Lambda^{0,3} \\
		\Lambda^{4}\left( T_x^* M \right) &=& \left( \Lambda^3\left( N_x^* \right) \otimes T_x^* \right) \oplus \left( \Lambda^2\left( N_x^* \right) \otimes \Lambda^2\left( T_x^* \right) \right) \oplus \left( N_x^* \otimes \Lambda^3\left( T_x^* \right) \right) \\
		&=& \Lambda^{3,1} \oplus \Lambda^{2,2} \oplus \Lambda^{1,3}.
	\end{eqnarray*}
	In this notation, the symbol \( \Lambda^{p,q}\) does \textit{not} refer to a decomposition with respect to a complex structure as it usually does. Instead, it refers to \( \Lambda^p(N^*_x) \oplus \Lambda^q(T_x^*)\) as written. We may write \( \rho_{x}\) in components as follows
	\begin{equation*}
		\rho_{x} = \rho^{3,0} + \rho^{2,1} + \rho^{1,2} + \rho^{0,3}.
	\end{equation*}
	The following observations are immediately apparent
	\begin{enumerate}
		\item \( \rho^{0,3} = 0\) since \( \iota^*\rho = 0\)
		\item \( \rho^{3,0} \neq 0\) since \( \iota^*\hat{\rho} > 0\).
	\end{enumerate}
	Since \( N_x\) and \( \Lambda^{2}\left( T_x^* \right)\) have the same dimension, it suffices to show that the kernel of the map \( n \mapsto \iota(n\lrcorner \rho)\) is trivial. Suppose that there exists a non-zero \( n \in N_x\) such that \( \iota^*(n \lrcorner \rho) = 0\). Then it must be the case that 
	\begin{equation*}
		n\lrcorner \rho \in \Lambda^{1,1} \oplus \Lambda^{2,0}.
	\end{equation*}
	Let \( J\) denote the almost complex structure on \( T_x M\) determined by \( \rho\). Then a formula from \cite{fosc17} tells us that 
	\begin{eqnarray*}
		*(n\lrcorner \rho) &=& -Jn \wedge \rho \in \Lambda^{2,2} \oplus \Lambda^{1,3} \\
		&\Rightarrow& Jn \wedge \rho^{3,0} = 0 \\
		&\Rightarrow& Jn \in N_x.
	\end{eqnarray*}
	Here, when we write \( -Jn \wedge \rho\) we mean the wedge product of the metric dual of \( -Jn\) with \( \rho\). This contradicts the fact that \( \rho + i\hat{\rho}\) must be a complex volume form since \( \rho\) is a stable 3-form. 
\end{proof}

Note that this lemma is not true if \( \iota^*\hat{\rho}\) vanishes at \( p\). If \( \iota P\) is a special Lagrangian submanifold in a Calabi-Yau manifold, then this condition is automatically satisfied. 

\begin{lemma}
	For all \( [\iota] \in \mathcal{C}_{\rho}\), the cokernel of \( dc_{[\iota]}:T_{[\iota]}\mathcal{S} \rightarrow \Omega^3(P) \) is isomorphic to \(\mathbb{R}\). 
\end{lemma}

\begin{proof}
	Note that 
	\begin{equation*}
		dc_{[\iota](n)} = \iota^*\left( \mathcal{L}_{n} \rho \right) = \iota^*\left( d(n\lrcorner \rho) \right)
	\end{equation*}
	since \( \rho\) is closed and where \( \mathcal{L}\) denotes the Lie derivative. So any 3-form in the image of \( dc_{[\iota]}\) is clearly exact. Furthermore \Cref{lem:isomorphism} says that the image of \( dc_{[\iota]}\) is \textit{all} exact 3-forms on \( P\). Therefore since any 3-form on \( P\) is closed, \( P\) being 3-dimensional, the result follows.
\end{proof}

Finally, we can state the Lagrange multipliers problem: 
\[ \text{Find the critical points of}\  \restr{f_{\tau}}{\tilde{C}_{\rho}}.\]
Next, we briefly review the finite-dimensional situation. 

\subsection{Finite dimensional review}

Suppose that \(\pi: E \rightarrow M\) is a rank \(k\) vector bundle over an \(n\)-dimensional manifold \(M\). Let \(s: M \rightarrow E\) be a section of \(E\) and \(f: M \rightarrow \mathbb{R}\) a function. Define \(Z = s^{-1}(0)\). Let \( ds:TM \rightarrow TE\) denote the full derivative of \(s\). Whenever we have a connection, let \( \pi_{V}\) denote the vertical projection and let \(Ds = \pi_V \circ ds\). If \(ds_p\) is constant rank \(r + n\) along \(Z\), so that \(ds_p\) has the same rank for each \(p \in Z\), then \(Z\) is a properly embedded submanifold of codimension \(r\) in \(M\).

\begin{theorem}[Lagrange multipliers theorem]
	In the context described above, a point \(p \in Z \subset M\) is a critical point of \(\restr{f}{Z}\) if and only if there exists \(\lambda \in E^*\) such that 
	\begin{equation*}
		df_p = \lambda \circ Ds_p.
		\label{eq:lme}
	\end{equation*}
	The set of \(\lambda\)s satisfying \Cref{eq:lme} is an affine space of dimension \((k - r)\) where \(k\) is the rank of \(E\) and \(r\) is the rank of \(Ds_p\) along \(Z\). 
	\label{thm:lmt}
\end{theorem}

\begin{definition} 
	The \textit{Lagrange function} is a function on \(E^*\) defined by: 
	\begin{equation*}
		\Lambda(\lambda) = f(\pi (\lambda)) - (\lambda \circ s \circ \pi) (\lambda).
	\end{equation*}
	\label{def:lagrange_function}
\end{definition}

For convenience, let \(S(\lambda) = (\lambda \circ s \circ \pi) (\lambda)\). Note that \(S\) is a section of \(E^{**}\). 

\begin{lemma}
	Let \(\Lambda\) be the Lagrange function. Then \(\lambda \in E^*\) is a critical point of \(\Lambda\) if and only if \(p = \pi(\lambda)\) is a critical point of \(\restr{f}{Z}\) and \(\lambda\) satisfies \(df_p = \lambda \circ Ds_p\). \end{lemma}

\begin{proof}
	We have
\begin{equation*}
	dS_{\lambda}(\dot{\lambda}) = (\pi_V (\dot{\lambda})) \circ s(\pi (\lambda)) + (\lambda \circ Ds \circ d\pi)(\dot{\lambda}).
	\label{eqn3}
\end{equation*}
Therefore, 
\begin{equation*}
	d\Lambda_{\lambda} (\dot{\lambda}) = (df_p \circ d\pi) (\dot{\lambda}) + dS_{\lambda}(\dot{\lambda}).
	\label{eqn4}
\end{equation*}
Next, since \(\lambda \in E^*\) satisfies \(df_p = \lambda \circ Ds_p\) we have \(s(\pi(\lambda)) = 0\). So \(d\Lambda\) becomes: 
\begin{equation*}
	d\Lambda_{\lambda} = (df_p + \lambda \circ Ds_p) \circ d\pi \equiv 0 \fall \dot{\lambda} \in T_pE^*
	\label{eqn5}
\end{equation*}
due to \Cref{thm:lmt}. 

Conversely, since \(d\Lambda_{\lambda} (\dot{\lambda}) = 0\) for all \(\dot{\lambda} \in T_{\lambda}E^*\), it is zero in particular for \(\dot{\lambda} \in \ker d \pi\). In this case, \(\pi_V(\dot{\lambda}) = \dot{\lambda}\) which means \(s(\pi(\lambda)) = 0\). Therefore \(df_p + \lambda \circ Ds_p = 0\) and \(p\) is a critical point of \(\restr{f}{Z}\) again by \Cref{thm:lmt}.

\end{proof}

Analogous results hold in infinite dimensions. See for example \cite{abmara12}.

\subsection{The perturbed SL equations}\label{perturbed-SL-equations}

In this section, we derive the Euler-Lagrange equations for the above Lagrange multipliers problem. These equations also appeared in \cite{dose09} but without much explanation. Returning to the context of \cref{lmp}, we define the Lagrange functional \(\Lambda \) by analogy with the finite-dimensional setting. Note that \( C^{\infty}(P)\) is dual to \( \Omega^{3}(P)\) in the sense that there is a pairing 
\begin{eqnarray*}
	C^{\infty}(P) \times \Omega^3(P) &\rightarrow& \mathbb{R} \\
	(f,\alpha) &\mapsto& \int_{P}f\alpha.
\end{eqnarray*}
Thus, the Lagrange functional \( \Lambda:C^{\infty}(P) \times \tilde{\mathcal{F}} \rightarrow \mathbb{R}\) is given by 
\begin{equation*}
	\Lambda(\lambda,\tilde{\iota}) = \int_{[0,1] \times P} \tilde{\iota}\hspace{.02cm}^* \tau + \int_{P} \lambda\iota^*\rho.
\end{equation*}
Here, \( \iota = \tilde{\iota}(1,\cdot)\). The tangent space of \( C^{\infty}(P) \times \mathcal{F}\) at \( (\lambda,\iota)\) is 
\begin{equation*}
	T_{(\lambda,\iota)}\left( C^{\infty}(P) \times \mathcal{F} \right) = C^{\infty}(P) \times \Gamma(\iota^*TM).
\end{equation*}
As with \( f_{\tau}\), the functional \( \Lambda\) is \( \tilde{\mathcal{G}}\)-invariant and therefore descends to a functional on \( \tilde{\mathcal{S}}\). Furthermore, its derivative \( d\Lambda\) is a \( \mathcal{G}\)-invariant one-form on \( C^{\infty}(P) \times \mathcal{F}\) given by 
\begin{equation*}
	d\Lambda_{(\lambda,\iota)}(l,n) = \int_{P}\iota^*(n\lrcorner \tau) + \int_{P}\lambda \iota^*(n\lrcorner d\rho) + \int_{P} l\iota^*\rho.
\end{equation*}
Note that since \( \rho\) is closed, Stokes' theorem gives
\begin{equation*}
	d\Lambda_{(\lambda,\iota)}(l,n) = \int_{P} \iota^*(n \lrcorner \tau)  \int_{P}d\lambda\wedge \iota^*\left( n\lrcorner \rho \right) + \int_{P}l\iota^*\rho.
\end{equation*}
If \( v \in \Gamma(T\iota P)\), then \( d\Lambda_{(\lambda,\iota)}(l,v) = d\Lambda_{(\lambda,\iota)}(l,0)\). Therefore \( d\Lambda\) descends to a one-form on \( C^{\infty}(P) \times \mathcal{S}\). The following notation was used in \cite{dose09} and will also be used in this paper.

\begin{definition}
	Let \( \alpha\) be a \( k\)-form on a manifold \( M\). Let \( \iota:P\rightarrow M\) be an embedding of a manifold \( P\) into \( M\) and suppose that \( \iota^* \alpha = 0\). Then \( \alpha\) defines an \( \iota^*T^*M\)-valued \( (k-1)\)-form on \( P\) called \( \alpha_{N}\) given by 
	\begin{equation*}
		\alpha_N(v_1, \dots, v_{k-1}) = \alpha\left( \iota_* v_1, \dots , \iota_* v_{k-1},\cdot \right)
	\end{equation*}
	for any \( (v_1, \dots, v_{k-1}) \in TP\). In the presence of a metric, \( \alpha_{N}\) takes values in \( N^* \iota\). 
	\label{def:aN}
\end{definition}

The above construction shows the following. 

\begin{proposition}
	We have \( d\Lambda_{(\lambda,\iota)}(l,n) = 0\) for all \( (l,n) \in T_{(\lambda,\iota)}\left( C^{\infty}(P) \times \mathcal{F} \right)\) if and only if the following two equations are satisfied
	\begin{align}
		\tau_{N} + d\lambda \wedge \rho_N &= 0 \label{eq:ele1}\\
		\iota^*\rho &= 0.
		\label{eq:ele}
	\end{align}
	\label{prop:ele}
\end{proposition}

\begin{definition}
	Let \( (M,\rho,\tau)\) be a 6-dimensional manifold equipped with a \( G_2\)-pair \( (\rho,\tau) \in \Omega^3(M) \times \Omega^4(M)\). Let \( P\) be an oriented, closed, 3-manifold. A pair \( (\lambda,\iota) \in C^{\infty}(P) \times \mathcal{F}\) which solves \cref{eq:ele1,eq:ele} is called a \textit{graphical associative} for reasons that will become clear in the next section.
	\label{def:grass}
\end{definition}

These are the \textit{Euler-Lagrange} equations for the Lagrange multipliers problem and will henceforth be referred to as the \textit{perturbed special Lagrangian (SL) equations}. The pair \( (\lambda,\iota)\) is a \textit{critical point} of the functional \( \Lambda\) if and only if \( (\lambda,\iota)\) satisfies these equations. By the Lagrange multipliers theorem, a pair \( (\lambda,\iota)\) is a critical point of \( \Lambda\) if and only if there exists \( \tilde{\iota} \in \tilde{\mathcal{F}}\) such that \( \iota = \tilde{\iota}(1,\cdot)\) and \( \tilde{\iota}\) is a critical point of \( \restr{f_{\tau}}{\tilde{C_{\rho}}}\). Note that the perturbed SL equations are \( \mathcal{G}\)-invariant. Thus the solutions may be thought of up to diffeomorphism. The following lemma appeared in \cite{hala82} and shows how the above setup relates to Calabi-Yau manifolds. 

\begin{lemma}[Harvey--Lawson]
	If \( (M,J,g,\Omega)\) is a 6-(real) dimensional Calabi-Yau manifold, then \( \im(\Omega)\) is a calibration whose calibrated submanifolds are called special Lagrangians. If \( \omega\) is the K\"ahler form for the Calabi-Yau metric \( g\), then a 3-submanifold \( L\) admits an orientation making it special Lagrangian if and only if \( L\) satisfies
	\begin{eqnarray*}
		\restr{\omega}{L} &=& 0 \\
		\restr{\V{Re}(\Omega)}{L} &=& 0.
	\end{eqnarray*}
	\label{lem:harvey-lawson-lemma}
\end{lemma}

Therefore it is easy to see that if, in the Calabi-Yau case, we let \( \tau = \frac{1}{2}\omega^{2}\) and \( \rho = \V{Re}(\Omega)\), a submanifold \( \iota P\) is special Lagrangian if and only if the embedding \( \iota\) together with \( \lambda = \V{constant}\) is a critical point of \( \Lambda\) (see \Cref{lem:grass}). Due to this fact we hope to define a Floer theory for the critical points of \( \Lambda\).  

However, we will need to consider solutions to \cref{eq:ele1,eq:ele} other than those for which \( \lambda\) is a constant because of the following theorem from \cite{mcle98} which shows that the critical points of \( \Lambda\) and therefore \( \restr{f_{\tau}}{\tilde{C}_{\rho}}\) cannot usually be isolated. 

\begin{theorem}[McLean]
	Suppose that \( M\) is a Calabi-Yau manifold and \( L \subset M\) is a special Lagrangian submanifold. Then the moduli space of nearby special Lagrangians is a smooth manifold of dimension equal to the first Betti number of \( L\). 
	\label{thm:mclean-theorem}
\end{theorem}

Of course, if we think of \( \restr{f_{\tau}}{\tilde{C}_{\rho}}\) as a functional on embeddings rather than on submanifolds (i.e. embeddings up to diffeomorphism) the critical points are not isolated since if \( \iota \) is a critical point, so is \( g^*\iota\) for any \( g \in \mathcal{G}\). In the following sections, we will prove that if one perturbs the special Lagrangian equations the critical points of \( \Lambda\) will become isolated (up to diffeomorphism). 

\section{Ellipticity and the \texorpdfstring{\( G_2\)}{g2}-cylinder}\label{ellipticity} 

In this section, we express the perturbed SL equations in terms of a section of an infinite-dimensional vector bundle and state the definition of the moduli space. Then we restrict our discussion to a slice of the action of the diffeomorphism group. We prove that the linearization of the relevant section is an elliptic operator by exploiting the relationship between 6 and 7 dimensions. Throughout, \( M\) will be a 6-manifold equipped with 

\begin{itemize}
	\item a smooth \( G_2\)-pair \( (\rho',\omega') \in \Omega^{3}(M) \times \Omega^2(M)\) such that both \( \rho'\) and \( \omega'\) are \textit{closed}
	\item a smooth \( G_2\)-pair \( (\rho,\tau) \in \Omega^{3}(M) \times \Omega^4(M)\) such that both \( \rho\) and \( \tau\) are \textit{closed and tamed by} \( (\rho',\omega')\).  
\end{itemize}

Also, let \( g\) denote the metric on \( M\) induced by the \( G_2\)-pair \( (\rho',\omega')\). More precisely, since \( \varphi'= \rho' + dt \wedge \omega'\) is a closed \( G_2\)-structure on \( \mathbb{R} \times M\), there is a corresponding metric on \( \mathbb{R} \times M\). This induces the metric, which we call \( g\), on \( \left\{ 0 \right\} \times M\) which we identify with \( M\). As before, \( P\) will be a closed 3-manifold. In this section, we deal with the space of smooth embeddings \( \iota:P \rightarrow M\) such that \( \iota\) belongs to a fixed homology class \( A\) and such that both \( \iota^* \hat{\rho}\) and \( \iota^*\rho'\) are positive. Let \( \mathcal{F}\) denote this space. 

Next, consider the vector bundle 
\begin{equation*}
	\mathcal{E} \rightarrow C^{\infty}(P) \times \mathcal{F}
\end{equation*}
whose fiber at \( (\lambda,\iota)\) is 
\begin{equation*}
	\mathcal{E}_{(\lambda,\iota)} = \Omega^{3}(P) \times \Omega^{3}(P,\iota^*T^* M).
\end{equation*}
Solutions to the perturbed SL equations then correspond to the zero set of the following section of this bundle. 
\begin{eqnarray*}
	L:C^{\infty}(P) \times \mathcal{F} &\rightarrow& \mathcal{E} \\
	(\lambda,\iota) &\mapsto& \left( \iota^*\rho,\tau_N + d\lambda \wedge \rho_N \right). 
\end{eqnarray*}

\begin{definition}
	Let \( \tilde{\mathcal{N}}\left( A,P;(\rho,\tau) \right)\) denote the \textit{moduli space} of solutions to the perturbed SL equations. That is, 
	\begin{equation*}
		\tilde{\mathcal{N}} = \tilde{\mathcal{N}}\left( A,P;(\rho,\tau) \right) = \left\{ (\lambda,\iota) \sth \iota \in \mathcal{F}, \iota^*\rho = 0, \tau_N + d\lambda \wedge \rho_{N} = 0 \right\}.
	\end{equation*}
	\label{def:moduli-space}
\end{definition}
Then \( \tilde{\mathcal{N}}\) can be identified with the zero set of \( L\) 
\begin{equation*}
	\tilde{\mathcal{N}} = \tilde{\mathcal{N}}\left( A,P;(\rho,\tau) \right) = L^{-1}(0).
\end{equation*}
As in the finite-dimensional case, we often want to neglect the Lagrange multiplier \( \lambda\). Let \( \pi_2:C^{\infty}(P) \times \mathcal{F} \rightarrow \mathcal{F}\) be the projection and let 
\begin{equation*}
	\tilde{\mathcal{M}} = \pi_2( \tilde{\mathcal{N}} ).
\end{equation*}
This is the \textit{moduli space of perturbed special Lagrangians}. A few remarks about this definition of the moduli space are in order. 

\begin{enumerate}
	\item The tilde above \( \mathcal{M}\) is there to remind us that we have not yet taken the quotient with \( \mathcal{G}\).
	\item At the moment, we do not need to keep track of the parameters \( A, P\) and \( (\rho,\tau)\) so we drop them from the notation. Later on, we will want to add them back in.
	\item Also at the moment, we have defined the moduli space in such a way that all its elements are smooth. In \cref{transversality} we will want to allow for non-smooth elements. The results of \cref{elliptic-regularity} will imply that if \( (\rho,\tau)\) are of class \( C^{\ell}\) then there exists a diffeomorphism \( \phi\) such that if \( (\lambda,\iota)\) is a \( C^{3}\) solution to the perturbed SL equations, \( (\lambda \circ \phi,\iota\circ \phi)\) is also of class \( C^{\ell}\). 
\end{enumerate}

\subsection{Restricting to a slice} The group of orientation-preserving diffeomorphisms of \( P\) acts freely on \( C^{\infty}(P) \times \mathcal{F}\) by composition (see \cite{cmmi91} for details) and \( \mathcal{S} = \mathcal{F}/\mathcal{G}\) is a Fr\'echet manifold. However, we will later need to work with Banach spaces and will want to identify elements in \( \left(C^{\infty}(P) \times \mathcal{F} \right)/ \mathcal{G}\) near a particular element with normal vector fields of varying regularity. Restricting to a slice of the action of \( \mathcal{G}\) removes this complication. 

The metric \( g\) corresponding to the \( G_2\) pair \( (\rho',\omega')\) induces an inner-product on the tangent spaces of \( C^{\infty}(P) \times \mathcal{F}\) as follows. Suppose that \( n_1, n_2 \in T_{\iota}\mathcal{F} = \Gamma(\iota^*TM)\) and \( l_1,l_2 \in T_{\lambda}C^{\infty}(P)\) . Define 
\begin{equation*}
	\langle n_1, n_2\rangle = \int_{P} g(n_1,n_2) \iota^* \rho'
	\label{eq:rho'-metric}
\end{equation*}
and 
\begin{equation*}
	\langle l_1,l_2\rangle = \int_{P} l_1 l_2 \iota^*\rho'.
\end{equation*}
These equations define a \( \mathcal{G}\)-invariant metric on \( C^{\infty}(P) \times \mathcal{F}\). 

\begin{definition}
	Using the exponential map corresponding to the metric \( g\) on \( M\), a neighborhood \( U_{(\lambda,\iota)}\) of \( (\lambda,\iota)\) can be identified with an open set in \( T_{(\lambda,\iota)} \left( C^{\infty}(P) \times \mathcal{F} \right)\). Under this identification, let \( S_{(\lambda,\iota)} \subset U_{(\lambda,\iota)}\) denote the subset corresponding to the orthogonal complement of the \( \mathcal{G}\)-orbit of \( (\lambda,\iota)\). We refer to \( S_{(\lambda,\iota)}\) as a \textit{local slice} for \( (\lambda,\iota)\) and it consists of pairs \( (l,n)\) where \( n\) is a normal (with respect to \( g\)) vector field along \( \iota P\).
	\label{def:slice}
\end{definition}

This definition makes sense because \( S_{(\lambda,\iota)}\) is transverse to the \( \mathcal{G}\)-orbit of \( (\lambda,\iota)\) in \( C^{\infty}(P) \times \mathcal{F}\). 

\begin{definition}
	Let \( \mathcal{N}^{S}\) denote the \textit{local moduli space of graphical associatives} given by
	\begin{equation*}
		\mathcal{N}^{S} = \mathcal{N}^{S}(A,P;(\rho,\tau)) = \left\{ (\lambda, \iota) \in L^{-1}(0) \sth (\lambda,\iota) \in S \right\}
	\end{equation*}
	where \( S\) is a local slice as described above. As before, also define
	\begin{equation*}
		\mathcal{M}^{S} = \pi_2(\mathcal{N}^S).
	\end{equation*}
	\label{def:slice-moduli-space}
\end{definition}

Let \( S\) be a local slice. Let \( \restr{L}{S}\) denote the restriction of the section \( L\) to \( S\).

\subsection{Ellipticity} Next, we prove 

\begin{theorem}
	The linearization of \( \restr{L}{S}\) at a point \( (\lambda, \iota) \in \mathcal{N}^{S}\) is a self-adjoint, elliptic operator.
	\label{thm:ellipticity}
\end{theorem}
\noindent First we review some definitions. 

\begin{definition}
	Let \( E,F\) be vector bundles over a manifold \( M\) and let \( T:\Gamma(E) \rightarrow \Gamma(F)\) be a differential operator of order \( k\). Let \((x,v) \in T^{*}M\) and \(e \in E_x\) be given. Find a smooth function \(g(x)\) and a section \(f \in \Gamma(E)\) such that \(dg_x = v\) and \(f(x) = e\). Then \textit{the principal symbol of \(T\)} is defined by 
	\begin{equation*}
		\sigma(T)(x,v)e = L\left( (g - g(x))^k f \right)(x) \in F_x
	\end{equation*}
	\label{def:symbol}
\end{definition}
%
\begin{definition}
	A symbol \( \sigma\) is called \textit{elliptic} if for all \( (x,v) \in T^* M\), the linear map \( \sigma(x,v): E_x \rightarrow F_x\) is an isomorphism. A differential operator \( T\) of order \( k\) is called \textit{elliptic} if its principal symbol is elliptic.
	\label{def:ellitpic-operator}
\end{definition}

Next, let \( (\lambda, \iota)\) be a graphical associative contained in a local slice \( S\). Let \( N \iota\) be the normal bundle of \( \iota P\) with respect to the metric determined by \( (\rho',\omega')\). Then the linearization of \( \restr{L}{S}\) is given by
	\begin{eqnarray*}
		D: C^{\infty}(P) \times \Gamma (N\iota) \rightarrow \Omega^3(P) \times \Omega^3 (P,  N^* \iota)\\
		D_{(\lambda, \iota)}(l,n) = \left( \iota^* d(n\lrcorner \rho) , d\left( n \lrcorner(\tau + d\lambda \wedge \rho) \right)_N + dl \wedge \rho_N \right).
	\end{eqnarray*}
In matrix form 
\begin{equation*}
	D_{(\lambda,\iota)}(l,n) = \begin{pmatrix}
		(d(\ \cdot \ )\wedge \rho)_N &&  d(\ \cdot\ \lrcorner (\tau+ d\lambda \wedge \rho))_N \\
		0 &&  \iota^* d(\ \cdot \  \lrcorner \rho)
	\end{pmatrix} 
	\begin{pmatrix}
		l \\
		n
	\end{pmatrix}.
	\label{hess(L)}
\end{equation*}

Note that \( D\) is an operator of order 1. Let \(\alpha \in T^*_x P\) and choose \((k,n) \in \mathbb{R} \times N_x \iota\). Let \(\sigma_{\tau,\rho}\) denote the symbol of \( D\). Furthermore, let \(\eta = \tau + d\lambda \wedge \rho\). Then \Cref{def:symbol} gives 
\begin{equation*}
	\sigma_{\tau,\rho}(x,\alpha) (k,n) = \left( \alpha \wedge(n \lrcorner \rho), (\alpha \wedge(n \lrcorner \eta))_N + k\alpha \wedge \rho_{N} \right)
	\label{rho tau symbol line}
\end{equation*}
or, in matrix form 
\begin{equation*}
	\sigma_{\tau,\rho}(x,\alpha) = \begin{pmatrix}
		\cdot \ (\alpha \wedge \rho_N) && \alpha \wedge \left( \cdot \ \lrcorner(\tau + d\lambda \wedge \rho) \right)_N \\
		0 && (\ \cdot \ \lrcorner \rho) \wedge \alpha
	\end{pmatrix} \begin{pmatrix}
		k \\
		n
	\end{pmatrix}.
	\label{eq:rho-tau-symbol-matrix}
\end{equation*}

In order to check that \(\sigma_{\rho,\tau}(x,\alpha)\) is an isomorphism for all choices of \((x,\alpha) \in T^* P\), we view our problem from the perspective of 7-dimensional \(G_2\) geometry since from this perspective, the symbol is significantly more simple. The next definition shows how the 6- and 7-dimensional perspectives are related: 

\begin{definition}
	Suppose that \((\lambda, \iota) \in C^{\infty}(P) \times \mathcal{F}\). Define 
	\begin{eqnarray*}
		\iota_{\lambda}: P &\rightarrow& \mathbb{R} \times M \\
		\iota_{\lambda}(p) &=& \left( \lambda(p) , \iota(p) \right).
	\end{eqnarray*}
	\label{def:iotalambda}
\end{definition} 
Then \(\iota_{\lambda}P\) is the graph of the function \(\lambda\) over \(\iota P\). The following lemma finally allows us to justify our term \textit{graphical associative}.

\begin{lemma}
	For any graphical associative \( (\lambda, \iota)\), \(\iota_{\lambda} P\) as defined above is an associative submanifold of \(X = \mathbb{R} \times M\) with respect to the \(G_2\) 4-form \(\psi = \tau + dt \wedge \rho\). 
	\label{lem:associative-graph}
\end{lemma}

\begin{proof}
	Note that any vector field tangent to \(\iota_{\lambda} P\) can be written in the form \(u = \overline{u} + d\lambda(\overline{u}) \partial_t\) where \(\overline{u}\) is a vector field on \(\iota P\). Let \((u,v,w)\) be three such vector fields tangent to \(\iota_{\lambda} P\) and calculate \(\iota_{\lambda}^*(u \lrcorner v \lrcorner w \lrcorner \psi)\). It will be easy to see that since \(\iota^*(\rho) = 0\) and \((\tau + d\lambda \wedge \rho)_N = 0\),  \(\iota_{\lambda}^*(u \lrcorner v \lrcorner w \lrcorner \psi)\) also vanishes. 
\end{proof}

The above lemma also allows us to prove that when \( M\) is an actual Calabi-Yau manifold, then the \textit{only} solutions to the perturbed SL equations are special Lagrangians. 

\begin{lemma}
	In the above situation, the \textit{only} solutions \( (\lambda,\iota)\) to the perturbed SL equations are those for which \( \iota P\) is a special Lagrangian submanifold of \( M\) and \( d\lambda = 0\). 
	\label{lem:grass}
\end{lemma}
\begin{proof}
	We already know that solutions with \( d\lambda = 0\) correspond to special Lagrangian submanifolds in \( M\), so it suffices to show that any solution \( (\lambda,\iota)\) satisfies \( d\lambda = 0\). However, note that if \( \iota_{\lambda} P\) is a graphical associative submanifold in \( \mathbb{R} \times M\), it must be volume-minimizing in its homology class. On the other hand, note that the volume of \( \iota_{0}P\) where \( \iota_{0}(p)= (0,\iota(p))\) is always less than or equal to the volume of \( \iota_{\lambda}\) where equality holds if and only if \( \lambda\) is constant. So the fact that \( \iota_{\lambda} P\) is associative means that \( \lambda\) must be constant.
\end{proof}

Just as special Lagrangian submanifolds are the critical points of a functional, associative submanifolds are also critical points of a functional. In fact, the situation is much simpler in this case since we do not have any constraints. Therefore the idea is to write the symbol of \(\sigma_{\tau,\rho}\) in terms of the symbol of the action functional for associative submanifolds which is easily seen to be elliptic. 

\subsection{The associative action functional}

A detailed description of the following material was given in \cite{sawa10}. Here it is slightly generalized. Suppose that \((X,\psi,\varphi')\) is a 7-dimensional manifold equipped with a closed, stable 3-form \( \varphi'\) and a closed \( \varphi'\)-tame 4-form \( \psi\). Let \( A \in H_{3}(X;\mathbb{R})\). Let \(P\) be a closed, oriented, 3-dimensional manifold. Then define 
\begin{equation*}
	\mathcal{F}(X) = \left\{ \iota: P \rightarrow X \sth \iota \in C^{\infty}\ \text{is an embedding,}\ \iota^*\varphi' > 0, \text{\ and\ } \iota \in A \right\}.
	\label{eq:emb7}
\end{equation*}
We also have the covering space versions of these objects exactly as in \cref{lmp}. The 4-form \(\psi\) defines a \(\mathcal{G}\)-invariant action functional on \(\tilde{\mathcal{F}}(X)\) by integration: 
\begin{equation*}
	f_{\psi}(\tilde{\iota}) = \int_{[0,1] \times P} \tilde{\iota}\hspace{.02cm}^* \psi.
	\label{eq:fpsi}
\end{equation*}
This descends to a \(\mathcal{G}\)-invariant, horizontal one-form on \(\mathcal{F}(X)\) 
\begin{equation*}
	(df_{\psi})_{\iota}(n) = \int_{P} \iota^*(n\lrcorner \psi).
	\label{eq:dfpsi}
\end{equation*}

Similar to the 6-dimensional case, \(\iota\) is a critical point of \(f_{\psi}\) if and only if \(\psi_N = 0\). In other words, if and only if \(\iota P\) is an associative submanifold of \(X\). Therefore, it's easy to see that the critical points of \(f_{\psi}\) are exactly the associative submanifolds of \(X\). As before, we can write down the symbol for the operator associated to the equation \(\psi_N = 0\) and restrict ourselves to a local slice of the action of the diffeomorphism group. 

Let \( N \iota\) denote the normal bundle of any embedding \( \iota:P \rightarrow X\) with respect to the metric determined by \( \varphi'\). Note that local slices of the action of \( \mathcal{G}\) on \( \mathcal{F}(X)\) can also be defined in this context.

Let \( (x,\alpha) \in T_x^* P\) and \( (x,n) \in N_x \iota \). Then 
\begin{equation*}
	\sigma_{\psi}(x,\alpha)(n) = \alpha \wedge (n \lrcorner \psi)_N.
	\label{eq:sigma-psi}
\end{equation*}
Now suppose that \(X = \mathbb{R} \times M\) where \(M\) is a 6-dimensional manifold with a closed, tamed, \(G_2\) pair. Let \((\lambda,\iota)\) be a critical point of \(\Lambda\). Then note that \(N_{x} \iota_{\lambda}  \cong \mathbb{R} \oplus N_x \iota \).

\begin{lemma}
	Under the identification \(N_{x} \iota_{\lambda} \cong \mathbb{R} \oplus N_x \iota \), we have
	\begin{equation*}
		\sigma_{\tau,\rho}(x,\alpha)(k, \overline{n}) = \sigma_{\psi}(x,\alpha)(k \partial_t+ \overline{n})
	\end{equation*}
	where \( (x,\alpha) \in T_x^*P\).
	\label{lem:6=7}
\end{lemma}

\begin{proof}
	Note that the image of \(\sigma_{\psi}(x,\alpha)\) lies in 
	\begin{eqnarray*}
		\Lambda^3(T^*_{x} P) \otimes N^*_{x}\iota_{\lambda} &\cong& \Lambda^3(T^*_{x}P) \otimes \left( \mathbb{R} \oplus N_x^* \iota\right)\\
		&\cong& \Lambda^3(T^*_{x}P) \oplus \left( \Lambda^3(T^*_{x} P) \otimes N^*_x \iota \right).
	\end{eqnarray*}
	So we can write \(\sigma_{\psi}\) as a matrix with respect to this splitting: 
	\begin{equation*}
		\sigma_{\psi}(x,\alpha)(k,\overline{n})= 
		\begin{pmatrix}
			A & B \\
			C & D
		\end{pmatrix}
		\begin{pmatrix}
			k \\
			\overline{n}
		\end{pmatrix}.
	\end{equation*}
	Let \( (x,\alpha) \in T_x^*P\) be arbitrary. When \(\overline{n} = 0\), we have 
	\begin{eqnarray*}
		\sigma_{\psi}(x,\alpha)(k, 0) &=& \alpha \wedge (k \partial_t \lrcorner \psi)_N \\
		&=& k \alpha \wedge \rho_N  \\
	\end{eqnarray*}
	which is an \(N^*_{x}\iota_{\lambda}\)-valued 3-form on \( P\). However, note that \(\rho\) is defined on \(\iota P\). So actually, we can view this as an \(N^*_x \iota \)-valued 3-form on \( P\). Therefore, we can conclude that 
	\begin{eqnarray*}
		A &=& \cdot \  \alpha \wedge \rho_{N} \\
		C &=& 0.
	\end{eqnarray*}
	Furthermore, when \(k = 0\), we have: 
	\begin{eqnarray*}
		\sigma_{\psi}(x,\alpha)(0,\overline{n}) = \alpha \wedge \left( \overline{n} \lrcorner \tau + \overline{n} \lrcorner(dt \wedge \rho) \right)_N.
	\end{eqnarray*}

	Since \( N^*_x \iota_{\lambda} \cong \mathbb{R} \oplus N^*_x \iota\), we can consider the projections onto the \( \mathbb{R}\) and \( N_x^* \iota\) factors. Similarly, we can consider the corresponding projections of \( \sigma_{\psi}(x,\alpha)(0,\overline{n})\). The first of these projections is \( B\) and the second is \( D\). We have:
	\begin{eqnarray*}
		B &=& \alpha \wedge (\cdot \ \lrcorner \tau + \cdot \ \lrcorner( d\lambda \wedge \rho))_N \\
		D &=& \alpha \wedge (\cdot \ \lrcorner \rho).
	\end{eqnarray*}
	After comparing this to equation \eqref{eq:rho-tau-symbol-matrix}, we see that the claim is proved.
	
\end{proof}

\begin{remark}
	Suppose that \(\iota P\) is an associative submanifold of a manifold \(X\) with \(G_2\)-structure. Let \(\left\{ e_1 , e_2, e_3 \right\}\) be an orthonormal basis for \(T_p \iota P\) and let \(\alpha = a_1 e^1 + a_2 e^2 + a_3 e^3\) be an arbitrary 1-form on \(\iota P\). Let \(n\) be an arbitrary vector normal to \(\iota P\) at \(p\). Then 
	\begin{eqnarray*}
		\left(e_3 \lrcorner e_2 \lrcorner e_3 \lrcorner \left( \alpha \wedge (n \lrcorner \psi) \right) \right)^{\flat} &=& a_1(e_1 \lrcorner e_2 \lrcorner n \lrcorner \psi)^{\flat} + a_2(e_3 \lrcorner e_1 \lrcorner n \lrcorner \psi)^{\flat} + a_3(e_2 \lrcorner e_1 \lrcorner n \lrcorner \psi)^{\flat} \\
		&=& a_1(e_3 \times n) + a_2(e_2 \times n) + a_3(e_1 \times n) \\
		&=& \alpha \times n
	\end{eqnarray*}
	where we have used the fact that \(e_1, e_2, e_3\) and \(n\) are all orthogonal, and the fact that \(\psi_N\) vanishes on \(\iota P\). In this sense, the symbol \(\sigma_{\psi}\) is just the dual of the symbol for the Dirac operator \(\dir = \sum_{i = 1}^{3} e_i \times \nabla_{e_i}\) whose kernel describes deformations of associative submanifolds and is elliptic, and self-adjoint.

	\label{elliptic7}
\end{remark}

\subsection{Self-adjointness}

The operator \( D_{(\lambda,\iota)}\) is self-adjoint in the following sense. Suppose that \( \iota: P \rightarrow \mathbb{R} \times M\) is an embedding. Let \( N\iota\) denote the normal bundle of this embedding and \( N^*\iota\) its dual. Then there is an isomorphism \( \Gamma(N^*\iota) \xrightarrow{\cong} \Omega^3(P,N^*\iota)\) given by the pairing:
\begin{equation*}
	\Gamma(N\iota) \otimes \Omega^3(P,N^*\iota) \rightarrow \mathbb{R} \qquad (n,\alpha) \mapsto \int_{P}n \lrcorner \alpha.
\end{equation*}
Let \( \iota_{\lambda}:P \rightarrow \mathbb{R} \times M\) be the embedding associated to a graphical associative \( (\lambda,\iota) \in C^{\infty} \times \mathcal{F}\) as in \Cref{def:iotalambda}. This identification, the above isomorphism, and the metric dual allow us to think of the operator \( D_{(\lambda,\iota)}\) as a map between \( N\iota_{\lambda}\) and itself. Explicitly, let \( n_1, n_2 \in \Gamma N \iota_{\lambda}\). Let \( \psi = \tau + dt \wedge \rho\) as usual. Then define 
 \begin{equation*}
	 D:\Gamma(N \iota_{\lambda}) \rightarrow \Gamma(N \iota_{\lambda}) \qquad n \mapsto \left( d(n\lrcorner \psi)_N \right)^{\flat}
 \end{equation*}
 where the metric dual is taken with respect to the metric \( g_{\psi}\) defined by \( \psi\). The following lemma was essentially proved by McLean in \cite{mcle98}, and again by Joyce in \cite{joyc16}. 

 \begin{lemma}
	 The operator \( D\) above is self-adjoint.
	 \label{lem:D-is-self-adjoint}
 \end{lemma}

Therefore, \Cref{lem:6=7} implies that \(\sigma_{\rho,\tau}\) is also elliptic and self-adjoint which is what we set out to prove. 

\section{Volume bounds}\label{volume-bounds}

In light of \Cref{lem:associative-graph}, we now have a natural way to identify solutions to the Lagrange multipliers problem with associative graphs in a cylinder \( \mathbb{R} \times M\). If \( \left( X, \psi, \varphi' \right)\) is a 7-manifold with a closed, tamed \( G_2\)-structure, and \( \iota:P \rightarrow X\) is a \( \psi\)-associative submanifold of \( X\), then \Cref{def:tamed-G2-structure} implies 
\begin{equation*}
	\V{Vol}_{g_{\psi}}(\iota P) \leq K \langle [\varphi'], [\iota]\rangle
\end{equation*}
for some positive constant \( K\), where \( [\iota]\) is the homology class represented by \( \iota\). We have an analogous volume bound for graphical associatives. 

\begin{lemma}
	Suppose that \( \left( M,\rho',\omega' \right)\) is a 6-dimensional manifold with a \( G_2\)-pair \( (\rho',\omega') \in \Omega^3(M) \times \Omega^2(M)\). Also suppose that \( (\rho,\tau)\) is a \( (\rho',\omega')\)-tame \( G_2\) pair and that \( (\lambda,\iota) \in \mathcal{N}(A,P;(\rho,\tau))\). Then both the volume of \( \iota P\) and \( \norm*{d\lambda}_{L^2}\) are bounded.
	\label{lem:cp-volume-bound}
\end{lemma}

\begin{proof}
	Suppose that \( (\lambda,\iota)\) is a graphical associative. Let \( \iota_{\lambda}\) denote the associated embedding \( \iota_{\lambda}:P \rightarrow \mathbb{R} \times M\), as before. We have the following metrics. 
	\begin{itemize}
		\item \( g_{M}\) : the metric on \( M\) given by the \( G_2\)-pair \( (\rho,\tau)\)
		\item \( g_{X}\) : the metric on \( X = \mathbb{R} \times M\) given by the \( G_2\)-structure \( \psi = \tau + dt \wedge \rho\)
		\item \( g = dt^2 + g_M\): the product metric on \( \mathbb{R} \times M\).
	\end{itemize}
	Since \( (\rho,\tau)\) is not necessarily an \( \V{SU}(3)\)-structure, \( g_X\) does not necessarily equal \( g\). Fix a point \( p \in P\) and let \( \left\{ e_i \right\}_{i=1}^{3}\) be a basis for \( T_{p} P\). Also let
	\begin{eqnarray*}
		\left\{ v_i \right\}_{i = 1}^{3} &=& \left\{ d\iota\left( e_i \right) \right\} \\
		\left\{ u_i \right\}_{i=1}^{3} &=& \left\{ d\iota_{\lambda}\left( e_i \right) \right\}.
	\end{eqnarray*}

	Choose \( \left\{ e_i \right\}\) so that \( v_{i}\) and \( \left\{ u_i \right\}\) are orthonormal with respect to \( g_M\) and \( g\) respectively. Denote the \( g_{M}\)-dual of \( v_{i}\) by \( v^{i}\). Similarly, denote the \( g\)-dual of \( u_i\) by \( u^{i}\). Let \( \mu\) be the function on \( \iota P\) defined by \( \lambda \circ \iota^{-1}\). Observe that

	\begin{equation*}
		u_i = v_i + d\mu(v_i)dt.
	\end{equation*}

	Let \( \V{vol}_{\iota_{\lambda}}\) denote the \( g\)-volume form on \( \iota_{\lambda} P\) and \( \V{vol}_{\iota}\) the \( g_M\)-volume form on \( \iota P\). Then we have \( \left( \V{vol}_{\iota_{\lambda}}\right)_{(t,x)} =  u^1 \wedge u^2 \wedge u^3\). We compute

	\begin{eqnarray*}
		u^1 \wedge u^2 \wedge u^3 &=& \left( v^1 + d\mu(v_1) dt \right) \wedge \left( v^2 + d\mu(v_2) dt \right) \wedge + \left( v^3 + d\mu(v_3) dt \right) \\
		&=& v^1 \wedge v^2 \wedge v^3 + d\mu(v_1) v^2 \wedge v^3 \wedge dt - d\mu(v_2)v^1 \wedge v^3 \wedge dt + d\mu(v_3)v^1 \wedge v^2 \wedge dt.
	\end{eqnarray*}

	Next, let \( L: \iota P \rightarrow \mathbb{R} \times M\) be the map defined by \( L(x) = \left( \mu(x),x \right)\). Then the above implies 

	\begin{equation*}
		L^*\left( \V{vol}_{\iota_{\lambda}} \right) = \V{vol}_{\iota} + *_{\iota}d\mu \wedge d\mu = \V{vol}_{\iota}\left( 1 + \abs*{d\mu}^2 \right)
	\end{equation*}

	Where \( *_{\iota}\) denotes the Hodge star on \( \iota P\) with respect to the metric induced by \( g_{M}\). Next, since \( \iota_{\lambda} P\) is associative by assumption, and since \( (\rho',\omega')\) tames \( (\rho,\tau)\), there exists a constant \( K > 0\) such that for all \( (t,x) \in \mathbb{R} \times M\),
	\begin{equation*}
		\left( \V{vol}_{\iota_{\lambda}} \right)_{(t,x)} \leq K\restr{\varphi'}{T_{(t,x)}\iota_{\lambda}P}.
	\end{equation*}

	Note that this is still true even though we are using the metric \( g\) instead of \( g_X\) to define the volume form \( \V{vol}_{\iota_{\lambda}}\). Since \( L^*(\varphi') = d\mu \wedge \omega' + \rho'\), we have
	\begin{eqnarray*}
		\int_{\iota P} \V{vol}_{\iota}\left( 1 + \abs*{d\mu}^2 \right) &\leq& K \int_{\iota P} d\mu \wedge \omega' + \rho' \\
		\Rightarrow \V{Vol}_{g_M}(\iota P) + \norm*{d\mu}_{L^2}^{2} &\leq& K \langle [d\mu \wedge \omega'],[\iota]\rangle + K \langle [\rho'],[\iota]\rangle \\
		&=& K\langle [\rho'],[\iota]\rangle.
	\end{eqnarray*}
	Clearly the \( L^2\)-norm of \( d\mu\) is bounded if and only if the \( L^2\)-norm of \( d\lambda\) is bounded.
\end{proof}



\section{Elliptic regularity}\label{elliptic-regularity}

In this section, we prove some regularity results that will be used in the proofs of the main theorems. The results in this section mainly follow from standard elliptic bootstrapping methods. We will work exclusively in the 7-dimensional setting for this section. In \cref{transversality} we will apply these to the 6-dimensional setting described above. 

\subsection{Set-up on \texorpdfstring{\( \mathbb{R}^{7}\)}{r7}} We begin by reviewing the so-called associator equation from \cite{hala82} which is a nonlinear PDE whose solutions are functions whose graphs are associative with respect to the standard \( G_2\)-structure on \( \mathbb{R}^{7}\). We then generalize this equation for non-flat \( G_2\)-structures. Let \( \mathbb{H}\) denote the quaternions and \( \mathbb{O}\) denote the octonions. Then identity 
\begin{equation*}
	\mathbb{R}^{7} \cong \V{Im}\mathbb{H} \times \mathbb{H} \cong \V{Im}\mathbb{O}.
\end{equation*} 
Furthermore, let \( \mathbf{1,i,j,k,e,ie,jk,ke}\) be the standard basis for \( \mathbb{O}\). In particular, any point \( x \in \mathbb{H}\) has the form 
\begin{equation*}
	x = x_0 + x_1 \mathbf{i} + x_2 \mathbf{j} + x_3 \mathbf{k}.
\end{equation*}
\begin{definition} 
	Suppose that \( x,y\) and \( z\) are octonions. Define the two- and three-fold \textit{cross products} to be 
	\begin{equation*}
		x \times y = -\frac{1}{2}\left( \bar{x}y - \bar{y}x \right)
	\end{equation*}
	and 
	\begin{equation*}
		x \times y \times z = \frac{1}{2}\left( x(\bar{y}z) - z(\bar{y}x) \right).
	\end{equation*}
	\label{def:cross-products}
\end{definition}

There are two important operators which are ingredients in the associator equation. 

\begin{definition}
	Let \( U\) be a domain in \( \V{Im}\mathbb{H}\) and let \( f:U \rightarrow \mathbb{H}\) be a \( C^{1}\) map. Then the \textit{Dirac operator} on \( f\) is defined to be 
	\begin{equation*}
		D(f) = -\frac{\partial f}{ \partial x_1}\mathbf{i} - \frac{\partial f}{\partial x_2}\mathbf{j} - \frac{\partial f}{\partial x_3} \mathbf{k}.
	\end{equation*}
	The \textit{first order Monge-Amp\'ere operator} on \( f\) is defined to be 
	\begin{equation*}
		\sigma(f) = \frac{\partial f}{\partial x_1} \times \frac{\partial f}{\partial x_2} \times \frac{\partial f}{ \partial x_3}.
	\end{equation*}
	\label{def:operators}
\end{definition}

\begin{theorem}[Harvey--Lawson]
	Let \( f:U \rightarrow \mathbb{H}\) be a \( C^1\) map. Then the graph of \( f\) is an associative submanifold of \( \V{Im}\mathbb{H} \times \mathbb{H}\) if and only if \( f\) satisfies 
	\begin{equation*}
		D(f) = \sigma(f).
	\end{equation*}
	\label{thm:associator-equation}
\end{theorem}

This equation is known as the \textit{associator equation} and is a first order, nonlinear PDE. It is well-known that \( C^1\) minimal submanifolds are smooth, and that associative submanifolds are minimal. Therefore a corollary to this theorem is that the graph of an associative, \( C^1\) map \( f:U \rightarrow \mathbb{H}\) is smooth. 

Let \( \psi_{0}\) denote the standard \( G_2\) 4-form on \( \mathbb{R}^{7}\). Let \( \beta \) be a 4-form on \( \mathbb{R}^{7}\) which is small enough so that \( \psi = \psi_0 + \beta\) is still stable. We now explain how to formulate a modified associator equation with respect to this new 4-form. To state this generalization we require the following constructions. Let 
\begin{itemize}
	\item \( g_0\) be the \( G_2\) metric on \( \mathbb{R}^{7}\) determined by \( \psi_0\) and
	\item \( g_{\psi}\) be the \( G_2\) metric on \( \mathbb{R}^{7}\) determined by \( \psi\). 
\end{itemize}

Also fix an orientation on \( \mathbb{R}^{7}\). To each stable 4-form, we also associate an \textit{associator} \( [\cdot,\cdot,\cdot]_{\psi}:\Lambda^3(\mathbb{R}^{7}) \rightarrow \mathbb{R}^7\) by requiring that 
\begin{equation*}
	g_{\psi}\left( [v_1,v_2,v_3]_{\psi},v_4 \right) = \psi(v_1,v_2,v_3,v_4) \quad \text{for all} \quad  v_1,v_2,v_3,v_4 \in \mathbb{R}^{7}.
\end{equation*}
Let \( [\cdot,\cdot,\cdot]_{\psi}^{0}:\Lambda^{3}(\mathbb{R}^{7}) \rightarrow \mathbb{R}^{7}\) denote the map determined by replacing \( g_{\psi}\) with \( g_0\) in the above formula. That is
\begin{equation*}
	g_{0}([v_1,v_2,v_3]_{\psi}^{0},v_4) = \psi(v_1,v_2,v_3,v_4) \quad \text{for all} \quad v_1,v_2,v_3,v_4 \in \mathbb{R}^{7}.
	\label{eq:modified-associator}
\end{equation*}
Then define \( \hat{\beta}(\cdot,\cdot,\cdot):\Lambda^3(\mathbb{R}^{7}) \rightarrow \mathbb{R}^{7}\) by the equation
\begin{eqnarray*}
	\psi(v_1,v_2,v_3,v_4) &=& \psi_{0}(v_1,v_2,v_3,v_4) + \beta(v_1,v_2,v_3,v_4) \\
	&=& g_0\left( [v_1,v_2,v_3]_{\psi_{0}},v_4 \right) + g_{0}\left( 2\hat{\beta}(v_1,v_2,v_3),v_4 \right).
\end{eqnarray*}
so that 
\begin{equation*}
	[\cdot,\cdot,\cdot]_{\psi}^{0} = [\cdot,\cdot,\cdot]_{\psi_{0}} + 2 \hat{\beta}(\cdot,\cdot,\cdot)
	\label{eq:associators}
\end{equation*}

Finally, thinking of \( \V{Im}\mathbb{O}\) as \( \V{Im}\mathbb{H} \oplus \mathbb{H}\), define projection operators \( \pi_{\mathbb{H}}\) and \( \pi_{\V{Im}{\mathbb{H}}}\).
\begin{eqnarray*}
	\pi_{\mathbb{H}}: \V{Im}\mathbb{O} &\rightarrow& \mathbb{H} \\
	\pi_{\V{Im}\mathbb{H}}: \V{Im}\mathbb{O} &\rightarrow& \V{Im}\mathbb{H}
\end{eqnarray*}
Then set 
\begin{eqnarray*}
	\hat{\beta}_{\mathbb{H}} &=& \pi_{\mathbb{H}} \circ \hat{\beta} \\
	\hat{\beta}_{\V{Im}\mathbb{H}} &=& \pi_{\V{Im}\mathbb{H}} \circ \hat{\beta}
\end{eqnarray*}

Next, suppose that \( f:U \rightarrow \mathbb{H}\) is a \( C^1\) map. Let \( f^i\) denote \( \frac{\partial f}{\partial x_i}\) for simplicity. Then the tangent space of the graph of \( f\) at each point is spanned by
\begin{eqnarray*}
	u &=& \mathbf{i} + f^1\mathbf{e} \\
	v &=& \mathbf{j} + f^2\mathbf{e} \\
	w &=& \mathbf{k} + f^3\mathbf{e} .
\end{eqnarray*}
Therefore we may think of \( \hat{\beta}_{\mathbb{H}}\) as a differential operator on \( f\) by defining 
\begin{equation*}
	\hat{\beta}_{\mathbb{H}}(f)\mathbf{e} = \hat{\beta}_{\mathbb{H}}(u,v,w).
\end{equation*}

\begin{theorem}
	As before, let \( U\) be a domain in \( \V{Im}\mathbb{H}\). Let \( f:U \rightarrow \mathbb{H}\) be a \( C^1\) map. Then the graph of \( f\) is associative with respect to \( \psi = \psi_0 + \beta\) if and only if it satisfies the equation 
	\begin{equation*}
		D(f) - \sigma(f) + \hat{\beta}_{\mathbb{H}}(f) = 0
	\end{equation*}
\end{theorem}

\begin{proof}
	After importing the above definitions, this proof proceeds similarly to the proof of the original theorem in \cite{hala82}. By definition, the graph of a function \( f\) is \( \psi\)-associative if and only if 
	\begin{equation*}
		[u,v,w]_{\psi} = 0
	\end{equation*}
	where \( u,v,w\) are vectors spanning the tangent space of the graph of \( f\). In this case, since 
	\begin{equation*}
		g_{\psi}([x_1,x_2,x_3]_{\psi},x_4) = \psi(x_1,x_2,x_3,x_4) = g_{0}([x_1,x_2,x_3]_{\psi}^{0},x_4)
	\end{equation*}
	for all \( v_1,v_2,v_3,v_4 \in \mathbb{R}^{7}\), the \( \psi\)-associative equation is equivalent to 
	\begin{equation*}
		[u,v,w]_{\psi}^{0} = [u,v,w]_{\psi_0} + 2 \hat{\beta}(u,v,w) = 0.
	\end{equation*}
	Also recall that \( \frac{1}{2}[x,y,z]_{\psi_{0}} = \V{Im}\left( x \times y \times z \right)\). Therefore, if the graph of \( f\) is \( \psi\)-associative, we have
	\begin{eqnarray*}
		\frac{1}{2}[u,v,w]^{0}_{\psi} &=& \V{Im}(u \times v \times w) + \hat{\beta}(u,v,w) \\
		&=& \V{Im}\left( \mathbf{i}\left( f^{2} \times f^{3} \right) + \mathbf{j}\left( f^3 \times f^1 \right) + \mathbf{k}\left( f^1 \times f^2 \right) \right) + \hat{\beta}_{\V{Im}\mathbb{H}}(f) \\
		&+& \left( \sigma(f) - D(f) \right)\mathbf{e} + \hat{\beta}_{\mathbb{H}}(f)\mathbf{e} = 0
	\end{eqnarray*}
	by formulas in \cite{hala82}. In particular, 
	\begin{equation*}
		\sigma(f) - D(f) + \hat{\beta}_{\mathbb{H}}(f) = 0
	\end{equation*}
	as desired. Conversely, suppose that \( \sigma(f) - D(f) + \hat{\beta}_{\mathbb{H}}(f) = 0\). Then 
	\begin{equation*}
		\frac{1}{2}[u,v,w]^{0}_{\psi} 
		= \V{Im}\left( \mathbf{i}\left( f^{2} \times f^{3} \right) + \mathbf{j}\left( f^3 \times f^1 \right) + \mathbf{k}\left( f^1 \times f^2 \right) \right) + \hat{\beta}_{\V{Im}\mathbb{H}}(f).
		\label{eq:eq1}
	\end{equation*}
	But 
	\begin{equation*}
		g_{0}([u,v,w]^0_{\psi},z) = g_{\psi}\left( [u,v,w]_{\psi},z \right) = 0
	\end{equation*}
	whenever \( z = u,v\) or \( w\) since \( [u,v,w]_{\psi}\) is orthogonal to \( u,v,\) and \( w\) with respect to the metric \( g_{\psi}\). Since 
	\begin{equation*}
		u = \mathbf{i} + f^1\mathbf{e},
	\end{equation*}
	this means that 
	\begin{equation*}
		g_0\left( [u,v,w]_{\psi}^{0},u \right) = g_{0}\left( [u,v,w]_{\psi}^{0},\mathbf{i} \right) + g_{0}\left( [u,v,w]_{\psi}^{0},f^{1}\mathbf{e} \right) = 0.
	\end{equation*}
	Similarly for \( v\) and \( w\). But since each term in \cref{eq:eq1} is contained in \( \V{Im}\mathbb{H}\) we conclude that
	\begin{equation*}
		g_0\left( [u,v,w]_{\psi}^{0},\mathbf{i} \right) = g_0\left( [u,v,w]_{\psi}^{0},\mathbf{j} \right) = g_0\left( [u,v,w]_{\psi}^{0},\mathbf{k} \right) = 0
	\end{equation*}
	Therefore
	\begin{equation*}
		[u,v,w]^{0}_{\psi} = 0
	\end{equation*}
	so the graph of \( f\) is associative as desired. 
\end{proof}

\subsection{Regularity and compactness}

Next, we apply elliptic bootstrapping methods to the modified associator equation to prove a regularity theorem and a compactness theorem for embeddings with bounded second fundamental form which will both play a crucial part in the transversality results to follow. Since the \( G_2\)-structures we are interested in in this paper are not necessarily torsion-free, their corresponding associative submanifolds are not necessarily minimal. This is why we cannot rely on regularity results about minimal submanifolds.

Recall that a subset \( S\) of a Riemannian manifold \( M\) is a \( C^{\ell}\) submanifold if it can be covered by \( C^{\ell}\) charts of \( M\). It is a standard result that a subset of a \( C^{\ell}\) manifold \( M\), where \( \ell \geq 1\), is a \( C^{\ell}\) submanifold if and only if it is the image of a \( C^{\ell}\) embedding. Furthermore, if \( S\) is the image of a \( C^1\) embedding \( \iota:P \rightarrow M\) and if \( S\) is also a \( C^{\ell}\) submanifold, then there exists a \( C^1\) diffeomorphism \( \phi\) of \( P\) such that \( \iota \circ \phi\) is a \( C^{\ell}\) embedding. See \cite{hirs12} for more details. 

\begin{theorem}[regularity]
	Let \( \ell \geq 2\) be an integer. Let \( (X,\psi)\) be a closed, oriented 7-manifold equipped with a stable \( C^{\ell}\) 4-form \( \psi\). Let \( \iota:P \rightarrow X\) be an embedding of a closed 3-manifold \( P\) into \( X\) so that \( \iota P\) is a \( C^{3}\) associative submanifold with respect to \( \psi\). Then \( \iota P\) is a \( C^{\ell}\) submanifold. In particular, if \( \psi\) is smooth, then so is \( \iota P\). 
	\label{thm:regularity}
\end{theorem}

\begin{proof}
	First note that at every point \( x \in \iota P\), there is a small neighborhood \( U\) of \( x\) in \( \iota P\) which can be regarded as a graph over \( T_x \iota P\). With this in mind, identify \( T_{x} \iota P\) with \( \V{Im} \mathbb{H}\), and \( N_x \iota\) with \( \mathbb{H}\). Let \( r\) be the distance from the origin in \( T_x X\). Let \( \alpha \in (0,1)\) be a number. It suffices to prove that a \( C^{2,\alpha}\) function \( f:U\rightarrow \mathbb{H}\) whose graph is associative with respect to a 4-form \( \psi = \psi_0 + \beta\) on \( \V{Im}\mathbb{H} \oplus \mathbb{H}\), where \( \beta\) is \( O(r)\) and  \( C^{\ell-1,\alpha}\), is in \( C^{\ell,\alpha}\). Assume that \( f\) is such a function. Then \( f\) satisfies 
	\begin{equation*}
		D(f) - \sigma(f) + \hat{\beta}_{\mathbb{H}}(f) = 0
	\end{equation*}
	due to the fact that \( \iota P\) is associative. Applying \( D\) to both sides of this equation yields 
	\begin{equation*}
		\Delta(f) - (D \circ \sigma)(f) + (D \circ \hat{\beta}_{\mathbb{H}})(f) = 0
	\end{equation*}
	where \( \Delta(f)\) is the Laplacian. Let \( L_1 = (D\circ \sigma)(f)\) which can be calculated explicitly 
	\begin{eqnarray*}
		L_1(f) &=& -\left( f^{11} \times f^2 \times f^3 \right)\mathbf{i} - \left( f^1 \times f^{12} \times f^3 \right) \mathbf{i} - \left( f^1\times f^{2} \times f^{13} \right) \mathbf{i} \\
		&-& \left( f^{21} \times f^2 \times f^3 \right)\mathbf{j} - \left( f^1 \times f^{22} \times f^3 \right) \mathbf{j} - \left( f^1 \times f^2 \times f^{23} \right) \mathbf{j} \\
		&-& \left( f^{31} \times f^2 \times f^3 \right) \mathbf{k} - \left( f^1 \times f^{32} \times f^3 \right) \mathbf{k} - \left( f^1 \times f^2 \times f^{33} \right) \mathbf{k}.
	\end{eqnarray*}

	We see that \( L_1\) is a nonlinear, second order operator with smooth coefficients. However, we can associate a \textit{linear} second order operator \( L_1^{f}\) to \( L_1\) with coefficients in \( C^{1,\alpha}\) \textit{depending on the first derivatives of \( f\)} simply by declaring the parts of \( L_1(f)\) that are first derivatives of \( f\) to be coefficients.

	Next, we address the operator \( (D \circ \hat{\beta}_{\mathbb{H}})\). It is useful to first get a better understanding of the operator \( \hat{\beta}_{\mathbb{H}}\). For the moment, let \( \left\{ e_i \right\}\) be a basis for \( \mathbb{R}^{7}\) and let \( \left\{ e^i \right\}\) be its dual basis. Then any \( C^{\ell-1,\alpha}\) 4-form on \( \mathbb{R}^{7}\) can be written as 
\begin{equation*}
	\beta = \sum A_{ijkl} e^{ijkl}
\end{equation*}
where \( A_{ijkl}:\mathbb{R}^{7} \rightarrow \mathbb{R}\) are \( C^{\ell-1,\alpha}\) functions. On the other hand, \( \hat{\beta}:\Lambda^3(\mathbb{R}^{7}) \rightarrow \mathbb{R}^{7}\) is given by 
\begin{equation*}
	\hat{\beta} = \sum A_{ijkl}e^{ijk}_{l}.
\end{equation*}

Next, let \( \mathbb{R}^{7} \cong \V{Im}\mathbb{H} \oplus \mathbb{H}\) as usual. Let \( x = (x_1,x_2,x_3)\) be coordinates on \( \V{Im}(\mathbb{H})\). 
At a point \( (x,f(x))\), we have 
\begin{equation*}
	\hat{\beta}(f)(x) = \sum A_{i}(x,f(x))B_{i}(x)e_i
\end{equation*}
where \( A_i\) come from the coefficients of \( \beta\) and \( B_i\) are products of the first partial derivatives of \( f\). We see that 
\begin{equation*}
	\hat{\beta}(f):\mathbb{R}^{3} \rightarrow \mathbb{R}^{7}.
\end{equation*}
Projecting onto the last 4 coordinates gives us the map 
\begin{eqnarray*}
	\hat{\beta}_{\mathbb{H}}(f):\mathbb{R}^{3} &\rightarrow& \mathbb{R}^{4} \\
	x &\mapsto& \sum_{i = 4}^{7}A_i(x,f(x))B_i(x)e_i.
\end{eqnarray*}
By the Leibniz rule, \( (D \circ \hat{\beta}_{\mathbb{H}})\) is a nonlinear differential operator with first and second order parts. Let \( L_2\) denote the second order part and \( L_3\) denote the first order part. Both \( L_2\) and \( L_3\) are nonlinear with coefficients in \( C^{\ell-2,\alpha}\). Similar to above, we can associate a \textit{linear} second order operator \( L_2^{f}\) to \( L_2\) whose coefficients depend on the first derivatives of \( f\). In this case however, the fact that the equation for \( \hat{\beta}_{\mathbb{H}}(f)\) involves the composition of \( f\), which is \( C^{2,\alpha}\) with the function \( A_{i}\), which is \( C^{\ell,\alpha}\), means that the linear, second order operator \( L_2^{f}\) has coefficients in \( C^{1,\alpha^2}\). The linear, first order operator \( L_3^{f}\) also has coefficients in \( C^{1,\alpha^2}\). 

Altogether, let 
	\begin{equation*}
		L^f = \Delta - L_1^f + L_2^f + L_3^f.
	\end{equation*}
	We must check that \( L^f\) is elliptic. The principal symbol of \( L^f\) is a 3 by 3 symmetric matrix. Its diagonal components are 
	\begin{eqnarray*}
		B_{11} &=& 1 + \left( \cdot \times f^1 \times f^2\right)\mathbf{i} + b_{11} \\
		B_{22} &=& 1 + \left( f^1 \times \cdot \times f^{3} \right) \mathbf{j} + b_{22} \\
		B_{33} &=& 1 + \left( f^1 \times f^2 \times \cdot \right)\mathbf{k} + b_{33}.
	\end{eqnarray*}
	Its off-diagonal components are 
	\begin{eqnarray*}
		B_{12} &=& \left( f^1 \times \cdot \times f^3 \right)\mathbf{i} + \left( \cdot \times f^2 \times f^3 \right) \mathbf{j} + b_{12} \\ 
		B_{13} &=& \left( f^1 \times f^2 \times \cdot \right)\mathbf{i} + \left( \cdot \times f^2 \times f^3 \right) \mathbf{k} + b_{13} \\
		B_{23} &=& \left( f^1 \times f^2 \times \cdot \right)\mathbf{j} + \left( f^1 \times \cdot \times f^3 \right) \mathbf{k} + b_{23}.
	\end{eqnarray*}

	The terms \( b_{ij}\) come from \( L_2^f\) and depend on the coefficients of \( \beta\), which are \( C^{\ell-1,\alpha^{2}}\) and \( O(r)\) and the first derivatives of \( f\) which can be made as small as we like since we can choose the neighborhood in the tangent plane of \( \iota P\) over which we are graphing to be as small as we like. Thus we see that the principal symbol of \( L^f\) is a small perturbation of the identity matrix and \( L^f\) is therefore elliptic with coefficients in \( C^{1,\alpha^{2}}\). 

	Standard results from Schauder theory and elliptic regularity allow us to conclude that since \( f\) solves \( L^f(f) = 0\), there is a number \( \alpha \in (0,1)\) for which \( f\) is a \( C^{3,\alpha}\) function on \( U\). That means that all the first derivatives of \( f\) are in fact contained in \( C^{2,\alpha}\) for some \( \alpha \in (0,1)\). So the same argument implies that \( f\) is in \( C^{4,\alpha}\) which means that \( L^f\) has coefficients in \( C^{3,\alpha}\). This process can be repeated up until the first derivatives of \( f\) are contained in \( C^{\ell-2,\alpha}\) and \( f\) is in \( C^{\ell-1,\alpha}\). Then \( f\) is in fact in \( C^{\ell}\). Further bootstrapping is prevented by the regularity of the coefficients of \( \beta\). 

\end{proof}

Next we prove a compactness result for embeddings with bounded second fundamental form. This will be invoked when applying the ``Taubes trick'' in \cref{transversality}. 

Suppose that \( X\) is a 7-manifold and let \( \iota: P \rightarrow X\) be an embedding. Fix any Riemannian metric on \( X\). For any \( x,y \in \iota P\) let \( d_X(x,y)\) denote the distance between \( x\) and \( y\) in \( X\). That is, the distance between \( x\) and \( y\) as points in \( X\). On the other hand, let \( d_{\iota P}(x,y)\) denote the distance between \( x\) and \( y\) in \( \iota P\). Let \( \V{II}(\iota) \) denote the second fundamental form of \( \iota P\) and \( \V{vol}(\iota P)\) denote the volume with respect to the metric. 

\begin{theorem}[compactness]\label{thm:compactness}
	Let \( \ell \geq 1\) be an integer, \( \alpha \in (0,1)\) a number, and let \( p \geq \frac{3}{1-\alpha}\). Let \( X\) be a closed 7-manifold equipped with sequence of stable 4-forms \( \psi_a\) converging to a stable 4-form \( \psi\) in the \( C^{\ell,\alpha}\) topology. Suppose also that \( \iota_{a}:P \rightarrow X\) is a sequence of \( W^{3,p}\) embeddings such that for all \( x,y \in P\), 
	\begin{enumerate}[label=(\roman*)]
	\item
		\begin{equation*}
			d_{X}\left( \iota_a (x),\iota_a(y) \right) \geq \frac{1}{K} d_{\iota_a P}\left( \iota_a(x),\iota_a(y) \right)
		\end{equation*}
		\item \begin{equation*}
				\norm*{\V{II}(\iota_a)}_{L^p} \leq K 
			\end{equation*}
		\item \begin{equation*}
				\V{vol}(\iota_a P) \leq K
			\end{equation*}
	\end{enumerate}
	for large enough \( a\) and some fixed constant \( K\). Also assume that \( \iota_a\) is \( \psi_a\) associative for each \( a\). Then there exists a subsequence \( \iota_b\) of \( \iota_a\) and a sequence of diffeomorphisms \( \phi_b\) of \( P\) such that \( \iota_b \circ \phi_b\) converges in the \( C^{\ell,\alpha}\)-topology to \( \iota\), an embedding also satisfying (i) and (ii). 
\end{theorem}

\begin{proof}
	The Sobolev embedding theorem guarantees that \( \iota_a \in C^{2,\alpha}\). Thus, \Cref{thm:regularity} implies that there is a sequence \( \phi_a\) of diffeomorphisms such that \( \iota_a \circ \phi_a \in C^{\ell+1,\alpha}\). In order to simplify notation, we therefore just assume that \( \iota_a \in C^{\ell +1,\alpha}\). The main theorem in \cite{breu12} implies that there exists an embedding \( \iota:P \rightarrow X\) satisfying conditions (i) and (ii), a subsequence \( \iota_{b}\) of \( \iota_a\), and a sequence of diffeomorphisms \( \phi_b\) such that \( \iota_b \circ \phi_b\) converges to \( \iota\) in the \( C^{2}\)-topology. The condition (i) above guarantees that \( \iota\) is also an embedding. Furthermore, \( \iota\) is \( \psi\)-associative and therefore also in \( C^{\ell+1}\). We just need to show that \( \iota_a\) converges in the \( C^{\ell,\alpha}\)-topology.

	Let \( x \in \iota P\) and let \( U \subset \iota P\) be an open set in \( \iota P\) containing \( x\) small enough so that it can be represented as a graph of a \( C^{\ell+1,\alpha}\) function \( f\). Let \( V = \iota^{-1}(U)\) and \( U_a = \iota_a(V)\). Then, since \( \iota_a\) converges to \( \iota\), each \( U_a\) can be represented as a graph of a function \( f_a\) over \( T_x \iota P\). Section 6 in \cite{breu12} together with the assumption (ii) guarantees that there is a uniform bound on the \( C^2\) norm of \( f_a\). Therefore there is a uniform bound on the \( C^{0,\alpha}\) norm of the first derivatives of \( f_a\). This in turn implies that there is a uniform bound on the \( C^{0,\alpha}\) norm of the coefficients of the operator \( L^{f_a}\). Thus, the Schauder estimates imply that there is a uniform bound on the \( C^{2,\alpha}\) norm of \( f_{a}\). This in turn implies that there is in fact a uniform bound on the \( C^{1,\alpha}\) norm of the coefficients of \( L^{f_a}\). This process can be repeated to obtain a uniform bound on the \( C^{\ell+1,\alpha}\) norm for \( f_a\). Therefore, after passing to a subsequence  \(f_a\) converges to \( f\) in the \( C^{\ell,\alpha}\) topology. Thus \( \iota_a\) also converges to \( \iota\) in the \( C^{\ell,\alpha}\)-topology as desired.

\end{proof}

\section{Transversality}\label{transversality}

In this section we prove that for a generic choice of a \( G_2\)-pair, the moduli space \( \mathcal{M}\) is a collection of isolated points. 

\subsection{Another slice} The group \( \mathbb{R}\) acts on \( C^{\infty}(P)/\mathcal{G}\) by addition of a constant. Let \( O_{\lambda}\) denote the \( \mathbb{R}\)-orbit of \( \lambda \in C^{\infty}(P)/\mathcal{G}\). Note that the \( L^{2}\)-orthogonal complement to the space of constant functions is 
\begin{equation*}
	I(P) = \left\{ h \in C^{\infty}(P)/\mathcal{G} \sth \int_{P} h\ \V{vol}_{P} = 0 \right\}.
\end{equation*}
The quotient space \( C^{\infty}(P)/ \sim\) where \( f \sim g\) if and only if \( f = g + c\) for some constant \( c\) can therefore be identified with \( I(P)\).

Let \( S\) be a local slice for the action of \( \mathcal{G}\) on \( C^{\infty}(P) \times \mathcal{F}\) as in \Cref{def:slice}. Let \( I \times S\) denote the set of pairs \( (\lambda,\iota)\) where \( \iota \in \mathcal{F} / \mathcal{G}\) and where \( \lambda \in C^{\infty}(P)/\mathcal{G}\) is \textit{also} a representative of an element in \( I(P)\), as above. Note that if \( (\lambda,\iota)\) is a graphical associative, so is \( (\lambda + c,\iota)\) where \( c\) is any constant. Therefore restricting to the slice \( I\) amounts to neglecting solutions up to translation in the \( G_2\) cylinder.

Let \( \mathcal{E}\) denote the bundle over \( I \times S\) whose fiber at \( (\lambda,\iota)\) is \( \Omega^{3}(P) \times \Omega^3(P,N^*\iota)\). Define 
\begin{align}
	\restr{L}{I\times S}:I\times S &\rightarrow \mathcal{E} \\ 
	(\lambda,\iota) &\mapsto \left( \iota^*\rho,\tau_{N} + d\lambda \wedge \rho_N \right).
	\label{eq:L}
\end{align}

\begin{definition}
	The \textit{local moduli space of perturbed special Lagrangians} is defined to be 
	\begin{equation*}
		\mathcal{M}^{I\times S} = \mathcal{M}^{I \times S}\left( A,P;(\rho,\tau)\right) = \restr{L}{I \times S}^{-1}(0).
	\end{equation*}
	\label{def:lmpsl}
\end{definition}

Let \( D_{(\lambda,\iota)}\) denote the linearization of \( \restr{L}{I\times S}\). We make the following observations about \( D_{(\lambda,\iota)}\).

\begin{observation}
	Note that before restricting to the slice \( I\), the operator \( D_{(\lambda,\iota)}\) always had a kernel of dimension \textit{at least} 1 since the set \( \left\{ (c,0) \right\}\) for constants \( c\) was always contained in the kernel. This set is no longer in the domain of \( D_{(\lambda,\iota)}\) after restricting the section \( L\) to the slice \( I\).  
\end{observation}

\begin{observation} 
	On the other hand, the kernel of \( D_{(\lambda,\iota)}\) at least sometimes contains other elements. Suppose that \( (M,\rho,\tau)\) is an actual Calabi-Yau manifold. Then \( (0,\iota) \in I \times S\) where \( \iota:P \rightarrow M\) is a special Lagrangian submanifold are solutions to the perturbed SL equations. In this case, elements \( (l,n)\) in the kernel of \( D_{(\lambda,\iota)}\) satisfy 

	\begin{enumerate}
	\item \( d \iota^*(n\lrcorner \rho) = 0\)
	\item \( d (n \lrcorner \tau)_N = 0\).
\end{enumerate}

Note that the first equation is true whenever \( n \lrcorner \rho\) is a closed 2-form on \( \iota P\). The second equation can be re-written using the formulas in \cite{fosc17} as 
\begin{eqnarray*}
	n \lrcorner \tau &=& * (n^{\sharp} \wedge \omega) \\
	&=& (Jn)^{\sharp} \wedge \omega.
\end{eqnarray*}
Then, \( d\left( (Jn)^{\sharp} \wedge \omega \right) = d(Jn)^{\sharp} \wedge \omega + (Jn)^{\sharp} \wedge d \omega\). This vanishes whenever \( (Jn)^{\sharp} \) is a closed 1-form since in the Calabi-Yau case, \( \omega\) is already closed. 
\end{observation}

These considerations motivate the following definitions. 

\begin{definition}
Let 
	\begin{equation*}
		\mathcal{R} = \left\{ (\rho,\tau) \in \mathcal{R}_{G_2} \subset \Omega^{3}(M) \times \Omega^{4}(M) \sth d\rho = 0, \quad d\tau = 0, \quad (\omega',\rho') \ \text{tames}\ (\rho,\tau)\right\}
	\end{equation*}
	be the \textit{space of parameters}.
	\label{def:space-of-parameters}
\end{definition}

\begin{definition}
	A pair \( (\rho,\tau) \in \mathcal{R}\) is called \textit{regular} (depending on \( A \) and \( P \)) if \(D_{(\lambda, \iota)}\) is surjective for all \( (\lambda, \iota) \in \mathcal{M}^{I,S}(A,P;(\rho,\tau))\). Let \( \mathcal{R}_{\V{reg}}\) denote the set of all regular pairs \( (\rho, \tau) \in \mathcal{R}\).
	\label{def:regular-parameters}
\end{definition}

The space \( \mathcal{R}\) should be compared to the space \( \mathcal{J}\) of \( \omega\)-tame almost complex structures in the symplectic setting. Note, that we do not require that \( \rho\) and \( \tau\) form an \( \V{SU}(3)\)-structure. We require that the forms \( \rho\) and \( \tau\) are closed so that the operator governing the deformation theory of the assocaitive submanifolds in \( \mathcal{R} \times M\) is self-adjoint and therefore has index zero. These conditions is also needed in order for the Lagrange multipliers problem described in \cref{lmp} to apply but are not needed to define the moduli space. 

The theorems in this section are consequences of the implicit function and Sard-Smale theorems which only apply in the Banach space setting. We address that next. 

\subsection{Banach space setup} Let \( W^{k,p}_{I \times S}\) denote the \( W^{k,p}\) completion of \( I \times S\). Similarly, let \( \mathcal{E}^{k,p}\) denote the \( W^{k,p}\) completion of the bundle over \( I \times S\) whose fiber at \( (\lambda,\iota)\) is \( \Omega^3(P) \times \Omega^3(P,N^*\iota)\). 

\begin{remark}
	Suppose that \( (\lambda,\iota) \in C^{\infty}(P) \times \mathcal{F}\) and \( S_{(\lambda,\iota)}\) is a local slice for \( (\lambda,\iota)\). Then since each element of \( S_{(\lambda,\iota)}\) corresponds to a normal vector field along \( \iota_{\lambda} P\), the \( W^{k,p}\) completion of \( S_{(\lambda,\iota)}\) consists of \( W^{k,p}\) normal vector fields along \( \iota_{\lambda} P\). Therefore if \( I\) is a local slice for \( \lambda \in C^{\infty}(P)/\mathcal{G} \), \( W^{k,p}_{I\times S}\) is identified with a subset of \( W^{k,p}(P,N\iota \times \mathbb{R})\).
\end{remark}

Let 
\begin{equation*}
	\mathcal{M}=\mathcal{M}\left( A,P;(\rho,\tau) \right) = \tilde{\mathcal{M}}/\sim
	\label{eq:moduli-space}
\end{equation*}
be the quotient space obtained by modding out by the action of \( \mathbb{R} \times \mathcal{G}\) on \( C^{\infty}(P) \times \mathcal{F}\) which is given by 
\begin{equation*}
	(c,\phi) \cdot (\lambda,\iota) = \left( (\lambda \circ \phi) + c, \iota \circ \phi \right).
\end{equation*}

This definition of \( \mathcal{M} \) corresponds to our previous definition, because after restricting to the slice \( I\), the projection from \( \mathcal{N}\) to \( \mathcal{M}\) is a bijection. Endow \( \mathcal{M}\) with the quotient subspace topology coming form the \( C^{\infty}\) topology on \( C^{\infty}(P) \times \mathcal{F}\).

The following is a consequence of the implicit function theorem for Banach spaces. 

\begin{theorem}
	The moduli space \( \mathcal{M}\left( A,P;(\rho,\tau) \right)\) is a set of isolated points whenever \( (\rho,\tau) \in \mathcal{R}_{\V{reg}}\).
	\label{thm:transversality-1}
\end{theorem}

\begin{proof}
	Fix a local slice \( I \times S\) for a particular element \( (\lambda_0,\iota_0) \in \mathcal{M}\). We first prove the result for the local moduli space \( \mathcal{M}^{I \times S}\). Let \( L:I \times S \rightarrow \mathcal{E}\) denote the map defined by \cref{eq:L}, dropping the restriction to \( I \times S\) from the notation. Since we want to apply the implicit function theorem, extend \( L\) to a map, also called \( L\):
	\begin{equation*}
		L: W^{k,p}_{I \times S} \rightarrow \mathcal{E}^{k-1,p}
	\end{equation*}
	where \( k \geq 3\) is an integer and \( p > 3\). If \( (\lambda,\iota) \in L^{-1}(0)\), then the Sobolev embedding theorem implies that \( (\lambda,\iota) \in C^{2,\alpha}\) so \Cref{thm:regularity} implies that \( (\lambda, \iota)\) is smooth. In this way, the local moduli space \( \mathcal{M}^{I\times S}\) is identified with the subset \( L^{-1}(0)\) in \( W^{k,p}_{I \times S}\). Next, note that any Riemannian metric on \( M\) induces a trivialization of \( \mathcal{E}^{k-1,p}\) over \( W^{k,p}_{I \times S}\) via parallel transport along geodesics in \( M\). Choose any such metric. For any \( (\lambda,\iota) \in W^{k,p}_{I \times S}\), let the parallel transport maps be denoted by
	\begin{equation*}
		T_0(\lambda,\iota):\mathcal{E}^{k-1,p}_{(\lambda_{0},\iota_0)} \rightarrow \mathcal{E}^{k-1,p}_{(\lambda,\iota)}
	\end{equation*}
	which is an isomorphism. Each fiber of \( \mathcal{E}^{k-1,p}\) is a Banach space. Therefore the map \( \bar{L}\) given by 
	\begin{eqnarray*}
		\bar{L}:W^{k,p}_{I \times S} &\rightarrow& \mathcal{E}^{k-1,p}_{(\lambda_0,\iota_0)} \\
		\bar{L}(\lambda,\iota) &=& T_{0}(\lambda,\iota)^{-1}\left( L(\lambda,\iota) \right)
	\end{eqnarray*}
	is a smooth map between Banach spaces. Also,
	\begin{equation*}
		\bar{L}(\lambda,\iota) = 0 \quad \iff \quad L(\lambda, \iota) = 0 \quad \iff \quad (\lambda,\iota) \in \mathcal{M}.
	\end{equation*}
	By the definition of \( \Sreg\), \( 0\) is a regular value of \( L\) and therefore also of \( \bar{L}\). So the implicit function theorem says that \( \mathcal{M}^{I \times S}\) is a Banach manifold with 
	\begin{equation*}
		T_{(\lambda,\iota)}\mathcal{M} = D\bar{L}(\lambda,\iota).
	\end{equation*}

	Since the operator \( D_{(\lambda,\iota)}\) is Fredholm with index zero whenever \( (\lambda,\iota) \in \mathcal{M}\), this shows that \( \mathcal{M}^{I\times S}\) is a set of points isolated in the \( W^{k,p}\)-topology. It remains to show that \( (\lambda,\iota)\) is isolated with respect to the \( C^{\infty}\) topology. Suppose that \( (\lambda,\iota) \in \mathcal{M}^{I \times S}\) is not isolated in the \( C^{\infty}\) topology. Then every \( C^{\infty}\) neighborhood of \( (\lambda,\iota)\) contains a second point \( (\lambda',\iota') \in \mathcal{M}^{I \times S}\). This is a contradiction since every \( W^{k,p}\) neighborhood of \( (\lambda,\iota)\) contains a \( C^{\infty}\) neighborhood. 

	Since the above is true for every local slice \( I \times S\), the result follows for \( \mathcal{M}\). 

\end{proof}

In the next theorem, we will want to work with \( C^{\ell,\alpha}\) parameters. Define 
\begin{definition}  
	\begin{equation*}
		\mathcal{R}^{\ell,\alpha} = \left\{ (\rho,\tau) \in C^{\ell,\alpha}\left( M,\Lambda^{3}T^*M \oplus \Lambda^4 T^*M \right) \sth d\rho = 0, \quad d\tau = 0, \quad (\omega',\rho')\quad \text{tames} \quad (\rho,\tau)\right\}
	\end{equation*}
	where \( \alpha \in (0,1)\). 
	\label{def:ell-alpha-parameters}
\end{definition}

\begin{definition} 
	Define the \textit{ local universal moduli space} \( \mathcal{M}^{I\times S}\left( A,P;\mathcal{R}^{\ell,\alpha} \right)\) to be the set of pairs \( \left( (\lambda,\iota),(\rho,\tau) \right)\) such that \( (\rho,\tau) \in \mathcal{R}^{\ell,\alpha}\) and \( (\lambda,\iota)\) is a solution to the perturbed SL equations contained in the \( C^{\ell,\alpha}\) completion of the local slice \( I \times S\).
	\label{def:local-universal-moduli-space}
\end{definition}

Next, we prove that \( \mathcal{M}^{I\times S}\left( A,P;\mathcal{R}^{\ell,\alpha} \right)\) is a Banach submanifold of \( W^{k,p}_{I \times S} \times \mathcal{R}^{\ell,\alpha}\). First, a Lemma.  

\begin{lemma}
	Let \( (\rho,\tau) \in \mathcal{R}^{\ell,\alpha}\) and suppose that \( (\lambda,\iota) \in C^{\ell,\alpha}\) solves the perturbed SL equations with respect to \( (\rho,\tau)\). Also assume that \( \iota^*\rho > 0\) as usual. Then for every nonzero \( (l,n) \in C^{\ell,\alpha}\left( P,\mathbb{R} \times N\iota \right)\), there exists a pair \( (\alpha,\beta)\) where \( \alpha\) is a closed 3-form on \( M\) and \( \beta\) is a closed 4-form on \( M\) such that 
	\begin{equation*}
		\int_{P}\iota(n\lrcorner \beta) + \int_{P} d\lambda \wedge \iota^*(n\lrcorner \alpha) + \int_{P} l\iota^*\alpha > 0.
	\end{equation*}
	\label{lem:linear-algebra-fact}
\end{lemma}

\begin{proof}
	\textbf{Case 1:} (\( n \neq 0 \) but \( l = 0\))
	In this case, we may choose \( \alpha = 0\) and show that the first term is nonzero, which follows from exactly the same arguments used to prove proposition A.2 of \cite{dowa17}. That is, choose \( p \in P\). Let \( U\) be a neighborhood of \( p\) and let \( V\) be a tubular neighborhood of \( U\). Since \( n \neq 0\), it must be nonzero somewhere in \( U\). Therefore, we can let \( f\) be a function supported in \( V\) such that \( df(n) \geq 0\) and \( df(n) > 0\) somewhere. Let \( \eta \) be a 3-form on \( M\) with \( \restr{\iota^*\eta}{U} = \restr{\V{vol}_{P}}{U}\) and \( n \lrcorner \restr{d\eta}{V} = 0\). Then, set 
	\begin{equation*}
		\beta = d(f\eta) = df \wedge \eta + f d\eta.
	\end{equation*}
	We have: 
	\begin{equation*}
		\int_P \iota^*(n \lrcorner (df \wedge \eta + f d\eta )) = \int_P df(n) \V{vol}_{P} > 0
	\end{equation*}
	as desired. 

	\textbf{Case 2:} (\( l \neq 0\) but \( n = 0\)). In this case, only the last term matters. But this case is easy since we can just choose \( U\) and \( V\) small enough so that \( l\) never vanishes. Then let \( f\) be a function on \( M\) such that \( \restr{\iota^*(f)}{U} = \restr{\frac{1}{l}}{U}\). Choose the 2-form \( \nu\) on \( M\) such that \( df \wedge \nu \) is \( \V{vol}_{P}\) when pulled back to \( U\). Then let \( \alpha = d(f \nu)\). 

	\textbf{Case 3:} (\( l \neq 0\) and \( n \neq 0\)). This is the same as in Case 1 since we can choose \( \alpha = 0\) and guarantee that the first term is positive. 
\end{proof}

\begin{proposition}
	Let \( I \times S\) be a local slice. Let \( \ell \geq 3\) and \( 3\leq k \leq \ell\) be integers. Let \( \alpha \in (0,1)\) be a number and \( p \geq \frac{3}{1-\alpha}\). Then \( \mathcal{M}^{I\times S}\left( A,P;\mathcal{R}^{\ell,\alpha} \right)\) is a separable Banach submanifold of \( W^{k,p}_{I \times S} \times \mathcal{R}^{\ell,\alpha}\).
	\label{prop:universal-moduli-space}
\end{proposition}

\begin{proof}
	Note that when \( 3 \leq k \leq \ell\), \( \mathcal{M}^{I\times S}\left( A,P;\mathcal{R}^{\ell,\alpha} \right)\) can be regarded as a subset of \( W^{k,p}_{I \times S} \times \mathcal{R}^{\ell,\alpha}\). For \( k \geq 1\), consider the (trivial) Banach space bundle 
	\begin{equation*}
		\mathcal{E}^{k-1,p} \rightarrow W^{k,p}_{I \times S} \times \mathcal{R}^{\ell,\alpha}
	\end{equation*}
	whose fiber over \( \left((\lambda,\iota),(\rho,\tau) \right)\) is 
	\begin{equation*}
		\mathcal{E}^{k-1,p}_{\left( (\lambda,\iota),(\rho,\tau) \right)} = W^{k-1,p}\left( P,\Lambda^3(T^*P)\oplus\left( \Lambda^3(T^*P)\otimes N^*\iota \right) \right)
	\end{equation*}
	Let \( L\) be a section given by
	\begin{eqnarray*}
		L:W^{k,p}_{I\times S}\times \mathcal{R}^{\ell,\alpha} &\rightarrow& \mathcal{E}^{k-1,p} \\
		\left( (\lambda,\iota),(\rho,\tau) \right) &\mapsto& \left( \iota^*\rho,\tau_N + d\lambda \wedge \rho_N \right).
	\end{eqnarray*}
	When \( k \geq 3\), \Cref{thm:regularity} applies, so \( \mathcal{M}^{I\times S}\left( A,P;\mathcal{R}^{\ell,\alpha} \right) = L^{-1}(0)\). In order to show that the universal moduli space is a submanifold, we need to show that the differential 
	\begin{align*}
	DL\left( (\lambda,\iota),(\rho,\tau) \right)&:W^{k,p}\left( P,\mathbb{R} \oplus N\iota \right) \times C^{\ell,\alpha}\left( M,\Lambda^3(T^*M) \oplus \Lambda^4(T^*M) \right) \\
	&\rightarrow W^{k-1,p}\left( P,\Lambda^3(T^*P)\oplus\left( \Lambda^3(T^*P)\otimes N^*\iota \right) \right)
	\end{align*}
	is surjective at every point \( \left( (\lambda,\iota),(\rho,\tau) \right) \in L^{-1}(0)\). The differential of \( L\) can be written in two components which are denoted by
	\begin{equation*}
		DL\left( (\lambda,\iota),(\rho,\tau) \right) = D_{(\lambda,\iota)} + D_{(\rho,\tau)}.
	\end{equation*}

	The same arguments from \cref{ellipticity} show that \( D_{(\lambda,\iota)}\) is a Fredholm operator. Therefore \( DL\left( (\lambda,\iota),(\rho,\tau) \right)\) has a closed image when \( \left( (\lambda,\iota),(\rho,\tau) \right) \in L^{-1}(0)\), so it suffices to show that its image is dense. 

	We first consider the case \( k =1\). Although in this case, \( L^{-1}(0)\) cannot be identified with the universal moduli space since the bootstrapping arguments of \cref{elliptic-regularity} don't apply, \( D_{(\lambda,\iota)}\) is still Fredholm so we can still prove that the operator \( DL\) is surjective on \( L^{-1}(0)\). Then, it will follow from elliptic regularity that \( DL\) is surjective for larger values of \( k\).

	Suppose that the image of \( DL\) at \( \left( (\lambda,\iota),(\rho,\tau) \right) \in L^{-1}(0)\) is \textit{not} dense. Then there exists a nonzero 
	\begin{equation*}
		\eta = (\eta_1,\eta_2) \in L^q\left( P,\Lambda^3(T^*P) \oplus \left( \Lambda^3(T^*P) \otimes N^* \iota \right) \right)
	\end{equation*}
	where \( \frac{1}{p} + \frac{1}{q} = 1\), such that 
	\begin{enumerate}[label=(\roman*)]
		\item \begin{equation*}
				\int_{P}\langle \eta,D_{(\lambda,\iota)}(l,n)\rangle\iota^*\rho' =0 \qquad \qquad  \fall (l,n) \in W^{1,p}\left( P,\mathbb{R} \oplus N\iota \right)
			\end{equation*}
		\item and 
			\begin{equation*}
				\int_{P} \langle \eta, D_{(\rho,\tau)}(\alpha,\beta)\rangle \iota^* \rho' = 0 \qquad \qquad \fall (\alpha,\beta) \in C^{\ell,\alpha}\left( M,\Lambda^3(T^*M) \oplus \Lambda^3(T^*M) \right)
			\end{equation*}
	\end{enumerate}

Now, (i) implies that \( \eta \in W^{1,p}\) since \( D_{(\lambda,\iota)}\) is Fredholm and self-adjoint. Furthermore, \( D^* \eta = 0\). In particular, the Sobolev embedding theorem implies that \( \eta\) is continuous. 

On the other hand, for any \( (\alpha,\beta) \in T_{(\rho,\tau)}\mathcal{R}^{\ell,\alpha}\), 
\begin{equation*}
	D_{(\rho,\tau)}(\alpha,\beta) = \left( \iota^*\alpha,\beta_N + d\lambda \wedge \alpha_N \right).
\end{equation*}
Therefore, letting \( l\) the metric dual of \( \eta_1\) on \( P\) and letting \( n\) be the normal vector field corresponding to \( \eta_2\), we have
\begin{equation*}
	\int_{P} \langle (\eta_1,\eta_2),D_{(\rho,\tau)}(\alpha,\beta)\rangle \iota^*\rho' = \int_{P}\iota^*(n\lrcorner \beta) + \int_{P} d\lambda \wedge \iota^*(n \lrcorner \alpha) + \int_{P}l\iota^*\alpha = 0.
\end{equation*}
by (ii). But \Cref{lem:linear-algebra-fact} shows that this is a contradiction. 

The fact that \( DL\left( (\lambda,\iota),(\rho,\tau) \right)\) is surjective whenever \( \left( (\lambda,\iota),(\rho,\tau) \right) \in L^{-1}(0)\) for \( k \geq 2\) follows from elliptic regularity. As mentioned before, as long as \( 3 \leq k \leq \ell\), the universal moduli space is equal to \( L^{-1}(0)\). Therefore the result follows from the implicit function theorem. Since \( W^{k,p}_{I\times S} \times \mathcal{R}^{\ell,\alpha}\) is separable, so is \( \mathcal{M}^{I \times S}\left( A,P;\mathcal{R}^{\ell,\alpha} \right)\).
\end{proof}

Lastly, we prove that the set of parameters for which the moduli space \( \mathcal{M}\left( A,P;(\rho,\tau) \right)\) is a zero dimensional manifold is ``generic''. More precisely, 

\begin{theorem}
	The set \( \mathcal{R}_{\V{reg}}\) is a residual subset of \( \mathcal{R}\).
\end{theorem}

\begin{proof}
	Let \( \mathcal{R}^{\ell,\alpha}_{\V{reg}}\) be the subset of \( \mathcal{R}^{\ell,\alpha}\) for which the associated operator \( D_{(\lambda,\iota)}\) from \Cref{prop:universal-moduli-space} is surjective whenever \( (\lambda,\iota)\) solves the perturbed SL equations. First, we show that \( \mathcal{R}^{\ell,\alpha}_{\V{reg}}\) is residual is \( \mathcal{R}^{\ell,\alpha}\). For that purpose, consider the projection 
	\begin{eqnarray*}
		\pi:\mathcal{M}^{I\times S}\left( A,P;\mathcal{R}^{\ell,\alpha} \right) &\rightarrow& \mathcal{R}^{\ell,\alpha} \\
		\left( (\lambda,\iota),(\rho,\tau) \right) &\mapsto& (\rho,\tau).
	\end{eqnarray*}
	By \Cref{prop:universal-moduli-space}, this is a map between separable Banach manifolds. Furthermore, 
	\begin{eqnarray*}
		d\pi\left( (\lambda,\iota),(\rho,\tau) \right):T_{\left( (\lambda,\iota),(\rho,\tau) \right)}\mathcal{M}^{I\times S}\left( A,P;\mathcal{R}^{\ell,\alpha} \right) &\rightarrow& T_{(\rho,\tau)}\mathcal{R}^{\ell,\alpha} \\
		\left( (l,n),(\alpha,\beta) \right). &\mapsto&(\alpha,\beta)
	\end{eqnarray*}
	Note that if \( \left( (l,n),(\alpha,\beta) \right)\) is in the tangent space of the local universal moduli space, then 
	\begin{equation*}
		D_{(\lambda,\iota)}(l,n) + D_{(\rho,\tau)}(\alpha,\beta) = 0.
	\end{equation*}
	Therefore both the kernel and cokernel of \( d\pi\left( (\lambda,\iota),(\rho,\tau) \right)\) are the same as that of \( D_{(\lambda,\iota)}\). So \( d\pi\left( (\lambda,\iota),(\rho,\tau) \right)\) is surjective if and only if \( D_{(\lambda,\iota)}\) is surjective. Therefore, \( \mathcal{R}^{\ell,\alpha}_{\V{reg}}\) is precisely the set of regular values of \( \pi\). The above holds for any local slice \( I \times S\). Thus, for large enough \( \ell\), the Sard-Smale theorem implies that \( \mathcal{R}^{\ell,\alpha}_{\V{reg}}\) is residual. 

	In order to extend the argument to smooth parameters, we must use the so-called Taubes trick, similar to what is done in order to prove that the set of regular \( \omega\)-tame almost complex structures is residual in the symplectic case. Let \( K > 0\) be a constant. Consider the set 
	\begin{equation*}
		\mathcal{R}_{\V{reg},K} \subset \mathcal{R}
	\end{equation*}
	of all smooth, stable, \( (\rho',\omega')\)-tame, \( G_2\) pairs \( (\rho,\tau)\) such that \( D_{(\lambda,\iota)}\) is surjective for every graphical associative \( (\lambda,\iota)\) which satisfies the following two conditions.

	\begin{enumerate}[label=(\roman*)]
		\item For any two points \( q,q' \in \iota P\) let \( d_{M}(q,q')\) denote the distance between them with respect to the metric induced by \( (\rho',\omega')\) on \( M\). Similarly, let \( d_{\iota P}(q, q')\) denote the distance between them with respect to the induced metric on \( \iota P\). Then the first condition for \( (\lambda,\iota)\) to be contained in \( \mathcal{R}_{\V{reg},K}\) is
			\begin{equation*}
				d_{M}(\iota(x),\iota(y)) \geq \frac{1}{K} d_{\iota P}(\iota(x),\iota(y))
			\end{equation*}
		\item Next, let \( p > 0\) be a number. Let \( \V{II}(\iota_{\lambda})\) denote the second fundamental form of the embedding \( \iota_{\lambda}:P \rightarrow \mathbb{R} \times M\) as in \Cref{def:iotalambda}. The second condition is that 
			\begin{equation*}
				\norm*{\V{II}(\iota_{\lambda})}_{L^p} \leq K.
			\end{equation*}
	\end{enumerate}

	Note that the conditions (i) and (ii) are both \( \mathcal{G}\)-invariant. Also note that since \( (\rho,\tau)\) is \( (\rho',\omega')\)-tame, the volume of every graphical associative is bounded by the same constant. Since we are no longer directly relying on Banach spaces, we no longer need to restrict ourselves to a slice of the action of \( \mathcal{G}\). Also, every graphical associative \( (\lambda,\iota)\) satisfies (i) and (ii) for some constant \( K > 0\). Therefore
	\begin{equation*}
		\mathcal{R}_{\V{reg}} = \bigcap_{K>0} \mathcal{R}_{\V{reg},K}.
	\end{equation*}
	Therefore we must prove that each \( \mathcal{R}_{\V{reg},K}\) is open and dense in the \( C^{\infty}\)-topology. First we prove that each is open by proving that its complement is closed. Assume that a sequence \( (\rho_a,\tau_a)\) converging to \( (\rho,\tau)\) in the \( C^{\infty}\) topology is contained in the complement of \( \mathcal{R}_{\V{reg},K}\). Then for each \( a\) there exists a \( (\rho_a,\tau_a)\) graphical associative \( (\lambda_a,\iota_a)\) satisfying (i) and (ii) such that \( D_{(\lambda_a,\iota_a)}\) is not surjective. Each graphical associative \( (\lambda,\iota)\) has an associated associative embedding \( \left( \iota_a \right)_{\lambda_a}\) and the conditions (i) and (ii) ensure that \Cref{thm:compactness} applies. 

	Thus, there exists a subsequence \( \left( \lambda_b,\iota_b \right)\) of \( \left( \lambda_a,\iota_a \right)\) and a sequence of diffeomorphisms \( \phi_b\) such that 
	\begin{equation*}
		\left( \iota_b \right)_{\lambda_b} \circ \phi_b \rightarrow \iota_\lambda
	\end{equation*}
	in the \( C^{\infty}\) topology. The condition (i) guarantees that not only is \( \iota\) an embedding, but \( \iota_{\lambda}\) is also \textit{graphical} since we used the distance in \( M\) instead of the distance in \( \mathbb{R} \times M\). The limit also satisfies both (i) and (ii). Furthermore \( D_{(\lambda,\iota)}\) is also \textit{not} surjective. Therefore \( (\rho,\tau)\) is not contained in \( \mathcal{R}_{\V{reg},K}\). Therefore \( \mathcal{R}_{\V{reg},K}\) is closed for any \( K > 0\). 

	Finally, we prove that \( \mathcal{R}_{\V{reg},K}\) is dense in \( \mathcal{R}\) with respect to the \( C^{\infty}\) topology. At this point, the argument is almost identical to the same argument in the symplectic case. See for example section 3.2 in \cite{mcsa04}. We include it here for completion. First, let \( \mathcal{R}^{\ell,\alpha}_{\V{reg},K}\) be the obvious \( C^{\ell,\alpha}\)-version of \( \mathcal{R}_{\V{reg},K}\). Note that 
	\begin{equation*}
		\mathcal{R}_{\V{reg},K} = \mathcal{R}^{\ell,\alpha} \cap \mathcal{R}.
	\end{equation*}
	Then \Cref{thm:compactness} still applies, so \( \mathcal{R}^{\ell,\alpha}_{\V{reg},K}\) is open in \( \mathcal{R}^{\ell,\alpha}\) with respect to the \( C^{\ell,\alpha}\) topology. Let \( (\rho_0,\tau_0) \in \mathcal{R}\). Since \( \mathcal{R}^{\ell,\alpha}_{\V{reg}}\) is dense in \( \mathcal{R}^{\ell,\alpha}\) by the Sard-Smale theorem, then there exists a sequence \( \left( \rho_{\ell},\tau_{\ell} \right) \in \mathcal{R}^{\ell,\alpha}_{\V{reg}}\)  (here the index \( \ell\) depends on \( (\ell,\alpha)\)) such that for all \( \ell \geq \ell_0\) for large enough \( \ell_0\), 
	\begin{equation*}
		\norm*{(\rho_0,\tau_0) - (\rho_{\ell},\tau_{\ell})}_{C^{\ell,\alpha}} \leq 2^{-\ell}.
	\end{equation*}
	Since \( (\rho_{\ell},\tau_{\ell}) \in \mathcal{R}^{\ell,\alpha}_{\V{reg},K}\) and since \(\mathcal{R}^{\ell,\alpha}_{\V{reg},K} \) is open in the \( C^{\ell,\alpha}\)-topology, there exists \( \epsilon_{\ell} > 0\) depending on \( (\ell,\alpha)\) such that for every \( (\rho,\tau) \in \mathcal{R}^{\ell,\alpha}\), 
	\begin{equation*}
		\norm*{(\rho,\tau) - (\rho_{\ell},\tau_{\ell})}_{C^{\ell,\alpha}} < \epsilon_{\ell} \qquad \Rightarrow \qquad (\rho,\tau) \in \mathcal{R}^{\ell,\alpha}_{\V{reg},K}.
	\end{equation*}
	Choose \( (\rho_{\ell}',\tau_{\ell}') \in \mathcal{R}\) to be any smooth element such that 
	\begin{equation*}
		\norm*{(\rho,\tau) - (\rho_{\ell},\tau_{\ell})}_{C^{\ell,\alpha}} < \min\left\{ \epsilon_{\ell},2^{-\ell} \right\}.
	\end{equation*}
	Then
	\begin{equation*}
		(\rho_{\ell}',\tau_{\ell}') \in \mathcal{R}^{\ell,\alpha}_{\V{reg},K} \cap \mathcal{R} = \mathcal{R}_{\V{reg},K}.
	\end{equation*}
	Therefore the sequence \( (\rho_{\ell}',\tau_{\ell}')\) converges to \( (\rho_0,\tau_0)\) in the \( C^{\infty}\) topology. Therefore we have shown that \( \mathcal{R}_{\V{reg}}\) is the intersection of a countable number of open, dense sets. Therefore it is residual, as desired. 

\end{proof}

\bibliographystyle{plain}

\bibliography{refs}

\end{document}